\documentclass[12pt,leqno]{article}
\usepackage[pdftex]{graphicx}
\usepackage{verbatim}

\newcounter{conjecture}
\newcounter{remark}
\newcounter{corollary}

\newenvironment{corollary}{\medskip{\bf Corollary \thecorollary.}
\addtocounter{corollary}{1}\em}{\rm}
\newcommand{\eqnsection}{
     \renewcommand{\theequation}{\thesection.\arabic{equation}}
     \makeatletter
     \csname @addtoreset\endcsname{equation}{section}
     \makeatother}
\newtheorem{theorem}{Theorem}
\newtheorem{lemma}{Lemma}
\newtheorem{proposition}{Proposition}

\newcommand{\lar}{\longrightarrow}

\newcommand{\aaa}{\alpha}

\newcommand{\lll}{\label}

\def \hh{\hskip1cm}
\def \ova{\overrightarrow}

\def \be{\begin{equation}}
\def \ee{\end{equation}}
\def \bt{\begin{theorem}}
\def \et{\end{theorem}}
\def \bea{\begin{eqnarray}}
\def \eea{\end{eqnarray}}
\def \bas{\begin{eqnarray*}}
\def \eas{\end{eqnarray*}}

\def \th{\theta}

\def \ang{\angle}
\def \tri{\triangle}

\def \vski{\vspace{12pt}}
\def \noi{\noindent}

\newcommand{\ttt}{\triangle}

\def \({\left(}
\def \){\right)}

\def \nn{\nonumber}

\def \bc{\begin{center} }
\def \ec{\end{center} }
\def \bs{\begin{slide} }
\def \es{\end{slide} }

\def\square{{\vcenter{\vbox{\hrule height.3pt
          \hbox{\vrule width.3pt height5pt \kern5pt
             \vrule width.3pt}
          \hrule height.3pt}}}}
\def\qed{{\hfill $\square$ \bigskip}}

\eqnsection
\begin{document}

\title{Newton Revisited: An excursion in Euclidean geometry}

\author{Greg Markowsky}

\maketitle

\begin{abstract}

This paper discusses the relationship between Kepler's Laws and Euclidean geometry. Many of the theorems are from {\it Principia} by Isaac Newton, but a more modern manner of presentation is adopted.

\end{abstract}

\section{Introduction}

The goal of this paper is to derive
Kepler's Laws of Planetary Motion from the Law of Universal
Gravitation using purely geometric methods. My motives for this are entirely aesthetic. There's
no question that modern calculus gets the job done, but there is a
certain thrill that comes from old-fashioned Euclidean geometry, with its
similar triangles and tangent lines. Except for a small number of instances, the reader is encouraged to forget all the calculus they know in order to better enjoy the mathematics. Since this paper deals with results that are known to be true with complete rigor, the style here will be quite informal and nonrigorous. As my advisor used to say, we're going to play fast and loose.

The idea for the paper came about by reading two books. The
identity of the first, {\it Principia} by Isaac Newton, should be
obvious. Newton's masterpiece was an inspiration to me when I
discovered it in graduate school, and is a feast for geometry lovers. Unfortunately, {\it
Principia} presents the modern reader with a few obstacles that get
in the way of the beautiful mathematics. To begin with, the style is
difficult and unfamiliar, and the translations that I have seen
contain words which are not in common use these days. Second,
Newton assumed many theorems about conics which apparently were well known to mathematicians of his day. Sadly,
studying the conics has fallen out of favor a bit in the time since
then, so that even professional mathematicians may have to do a bit of work on their own to make it through {\it Principia}, as I did.

Given these difficulties, it is natural to attempt to "translate"
Newton's work into modern notation, with background material
supplied where necessary. The most notable recent attempt at this
that I know of was by Richard Feynman, and presented in a lecture to
students at the California Institute of Technology(in fact, the core of the argument that Feynman used had already been discovered by the great James Maxwell, and was published in \cite{maxwell}, though Feynman was probably unaware of it). Feynman's lecture has survived in the form of \cite{feynman}.
This is a very enjoyable book, and
Feynman's discussion is ingenious, but it fell short of real
satisfaction for me for two reasons. First of all, Feynman does
not derive Kepler's Third Law in its entirety from Newton's laws. Perhaps his methods could lead to a derivation, but it isn't mentioned in the book. Secondly, and more seriously, Feynman went out of his way to avoid dealing too much
with the geometrical properties of the conics. Certainly a person
has the right to dislike the conics if they choose, but I have
always ascribed more to the Archimedean school of thought which
contends that the greatest joy in physics lies in the
wonderful geometrical problems that arise. In other words, the
motion of the planets around the sun gives us a great excuse to
study the conics.

This paper, then, is a record of my attempt to understand Newton's work on planetary orbits. My argument departs from his at some point, but up to that point the paper is largely just a retelling of selected pieces of {\it Principia}. The next section is introductory material on conics. The third section follows Newton's work, with auxiliary lemmas added where needed. The first theorem in the fourth section is also from {\it Principia}, while the rest of the section gives the way I came up with to deduce Kepler's first and third laws. This part of the paper can be considered original, to my knowledge. The fifth section gives a solution to a natural problem in mechanics which is found to be extremely simple using the methods in the fourth section. If the reader can obtain from this paper $10\%$ of the enjoyment that I felt while studying {\it Principia}, they can consider it time well spent.

\section{Conics}

It is expected that the reader has some familiarity with the conics. For the sake of completeness, however, I have included essentially everything relevant below. A reader with a good working knowledge of the conics can safely skim this section. Before we look at the conics, I should mention some possibly nonstandard bits of terminology that I will use. The {\it measure} of an arc on a circle is defined to be the magnitude of the
angle it subtends at the center of the circle. For example, the measure of arc $AB$(abbreviated as $m(AB)$) below is
$35^\circ$.

\includegraphics[width=100mm, height=70mm]{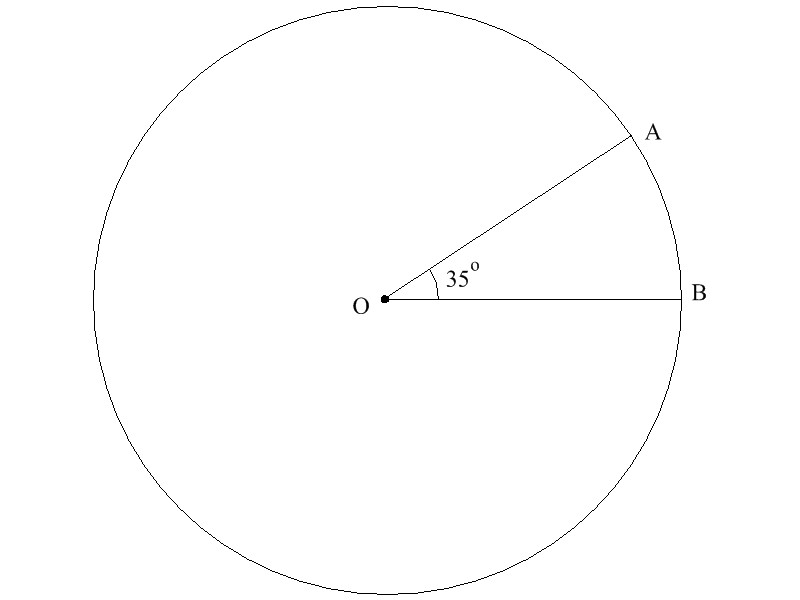}

Furthermore, when
an angle inside a circle subtends a pair of arcs in both directions,
we say that the angle {\it covers} the arcs. For example, in the
picture below, angle $\aaa$ covers arcs $AB$ and $CD$.

\includegraphics[width=100mm, height=70mm]{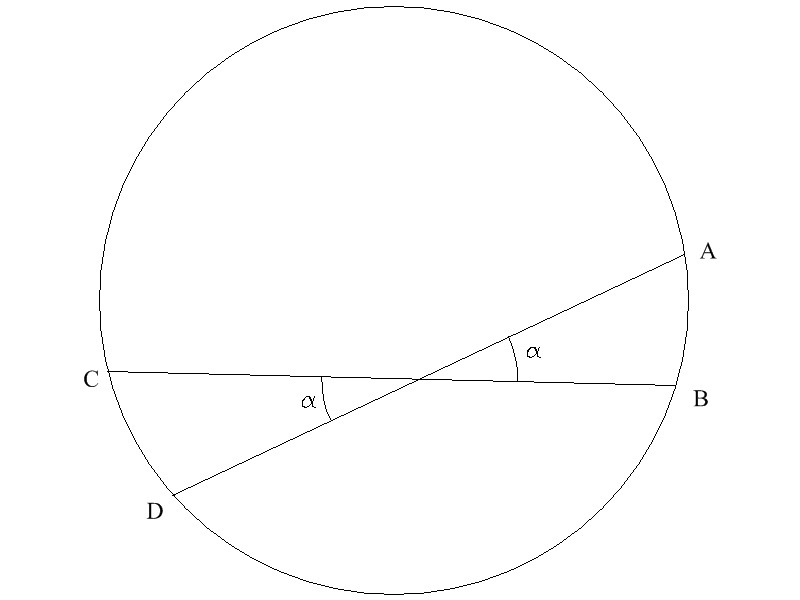}

For reasons which I am at a loss to explain, the following
theorem is rarely given in full generality when it is presented in high schools.

\begin{theorem} Suppose that two lines intersect at an angle $\aaa$ either inside or on a circle.
Let $AB$ and $CD$ be the arcs covered by $\aaa$. Then

\be \lll{}m(AB) + m(CD) = 2\aaa \ee

\end{theorem}

{\bf Proof:} Draw the lines parallel to the original two lines but
which pass through the center of the circle $O$. Let $A', B', C',
D'$ be the points of intersection of these new lines with the
circle, as shown below.

\includegraphics[width=100mm, height=70mm]{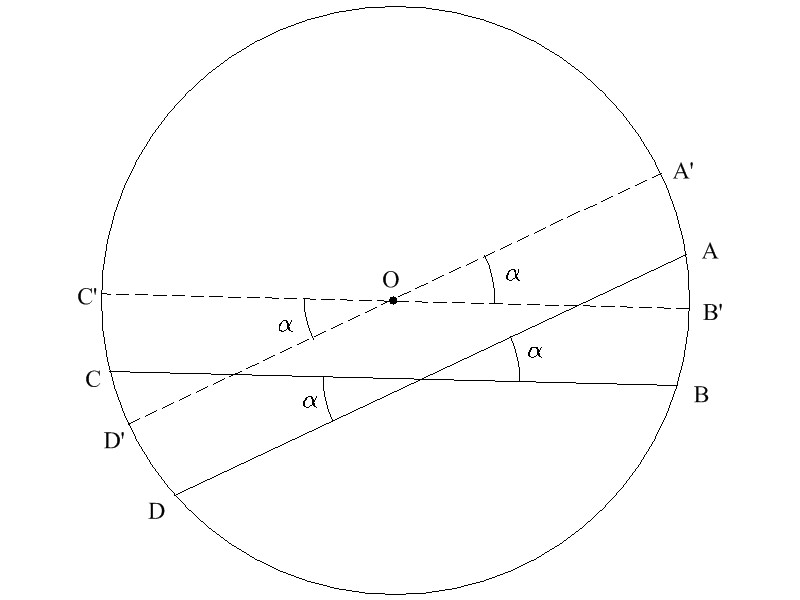}

We have

\be \lll{} m(AB) + m(CD) = m(AB') + m(C'D) = m(A'B') + m(C'D') \ee

As $m(A'B') + m(C'D') = 2\aaa$ by definition, we are
done. \qed

Note that this proof works just as well if the angle lies on the
circle, and includes the case where one of the lines is a tangent to
the circle.

The tangent to a curve more general than a circle must be defined. Let $O$ and $P$ be two points on a curve which are close to each other, and let $P$ be fixed. Draw the line
containing both $O$ and $P$. As we let $O$ approach $P$, if the line
containing $O$ and $P$ gets closer and closer to a fixed line $T$,
we define $T$ to be the tangent at $P$.

\includegraphics[width=100mm, height=70mm]{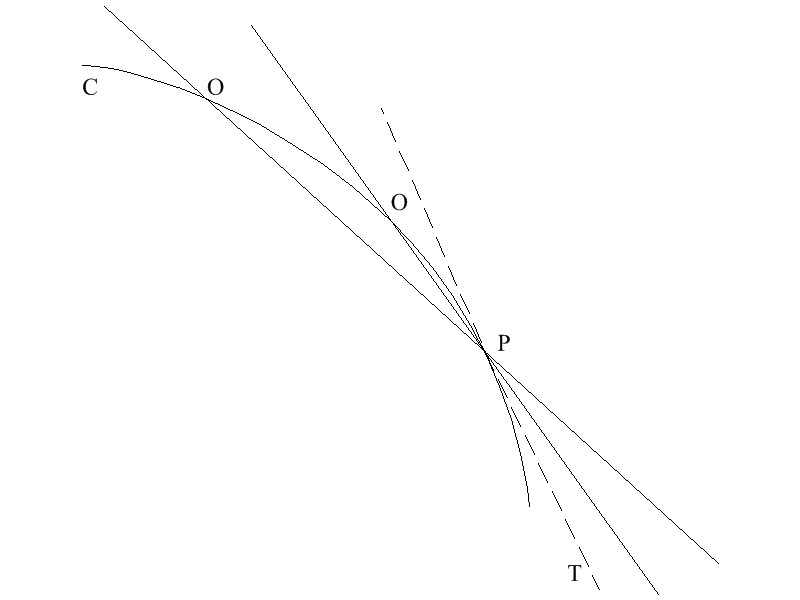}

\noi We aren't
going to worry about whether such a $T$ exists, we'll just assume it
does(as it does) for the curves we care about.

Finally, I will write things like
"the slope of the cone" and "the slope of the plane". These are not
standard terminology, but should cause no confusion. Given a cone,
orient it so that the vertex is pointing directly up. The vertical
line through the vertex is the axis of the cone, and we
define the slope of the cone to be the ratio $\frac{OV}{OA}$ below.
Note that $O$ is the center of the base, and $A$ is a point on the
outside of the base.

\includegraphics[width=100mm, height=70mm]{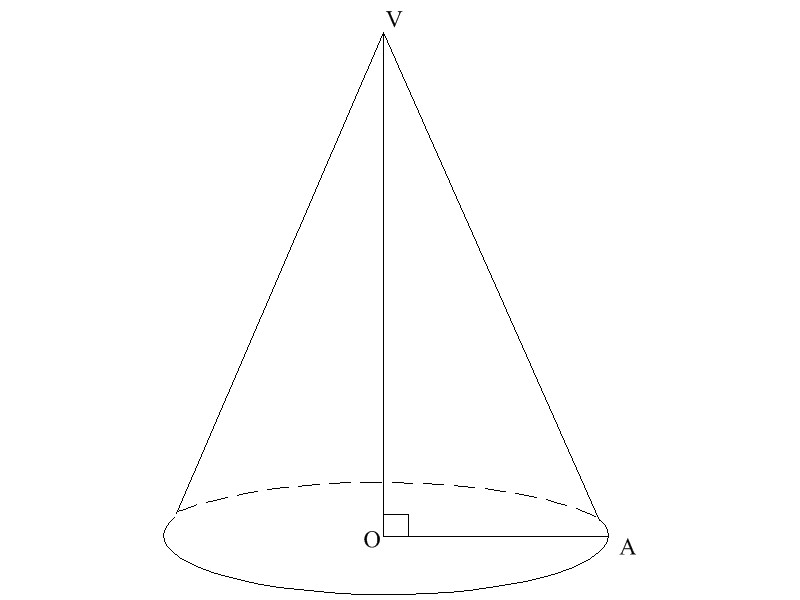}

To define the slope of a plane $P$, begin by drawing a vertical line
$L$ passing through the plane at a point $P$. Choose a point $O$ on
$L$ but not on the plane, and let $A$ be the closest point on the
plane so that $OA$ is perpendicular to $L$. Then the slope of the
plane is defined to be $\frac{OP}{OA}$.

\includegraphics[width=100mm, height=70mm]{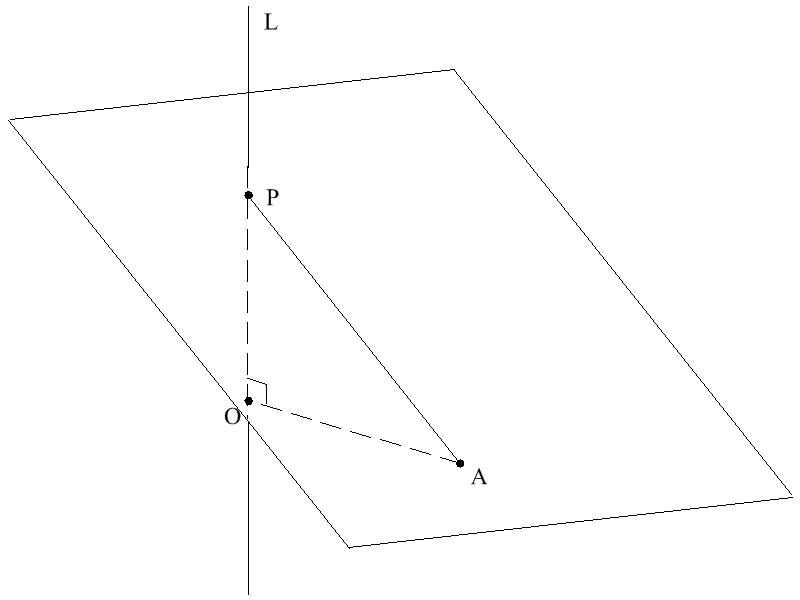}

The following lemma will help in dealing with tangents.

\begin{lemma} \lll{tangent} Suppose that a circle is tangent to a pair of lines
with points of tangency $a$ and $b$, and that the pair of lines meet
at $n$. Let $o$ be the center of the circle. Then
$\ova{ao}+\ova{bo}$ is parallel to $no$.
\end{lemma}
\includegraphics[width=100mm, height=70mm]{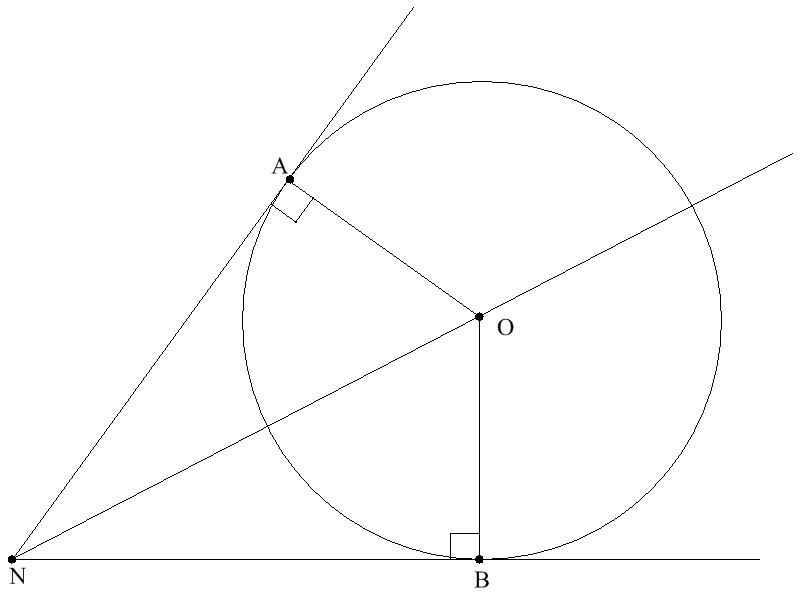}

{\bf Proof:} The proof is immediate, since the entire picture is
symmetric around the angle bisector $no$. \qed

Suppose we have two circles of different radii with the smaller
contained in the larger.

\includegraphics[width=100mm, height=70mm]{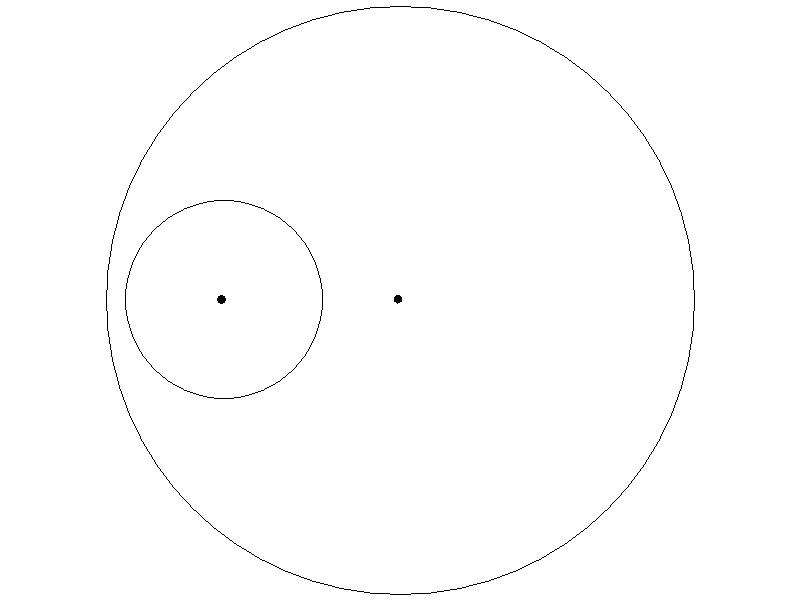}

We are allowing in this setup that the smaller is internally
tangent at a point to the larger. Let $E$ be the set of all points
which are centers of circles tangent to both of our original
circles. Then $E$ is a curve, and this curve is an {\it
ellipse}.

\includegraphics[width=100mm, height=70mm]{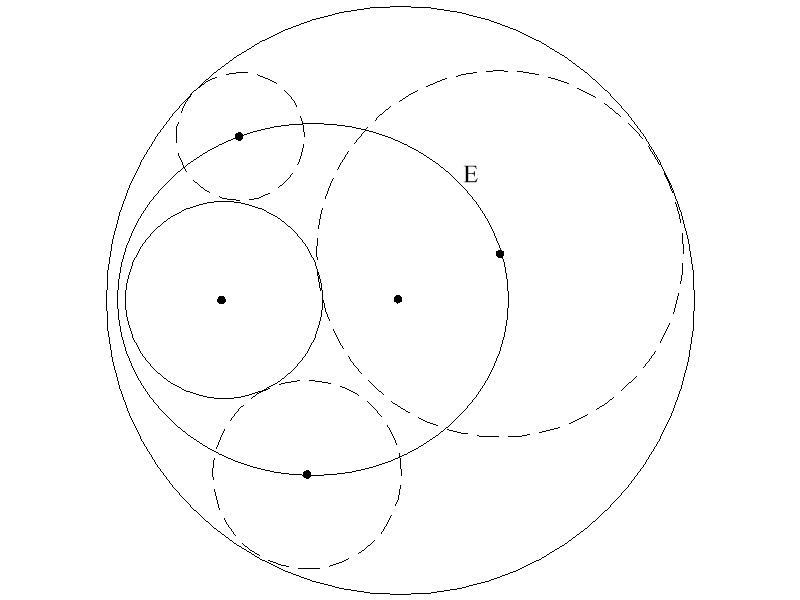}

The two centers of the original circles are known as the {\it foci}
of the ellipse. The following proposition gives the property of
ellipses which is usually used to characterize them.

\begin{proposition} \lll{sum} The sum of the distances from the foci of an
ellipse to any point on the ellipse is equal to the sum of the radii
of the two original circles.
\end{proposition}
{\bf Proof:} Examine the picture below.

\includegraphics[width=100mm, height=70mm]{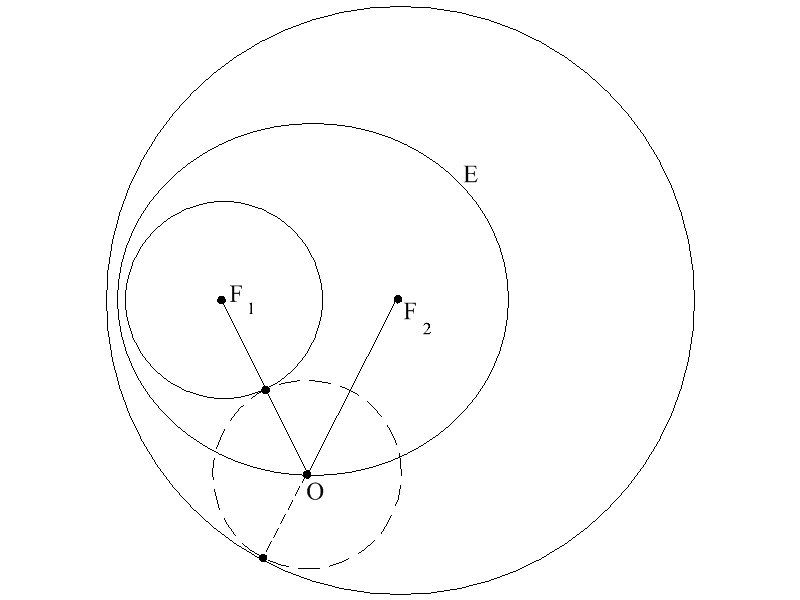}

$F_1$ and $F_2$ are the centers of the two fixed circles of radius
$r_1$ and $r_2$, and $O$ is a point on the ellipse whose
corresponding circle has radius $r_3$. We see that $F_2 O = r_2 -
r_3$, and $F_1 O = r_1 + r_3$. Thus, $F_2 O + F_1 O = r_1 + r_2$, a
constant. \qed

And now we have the all-important reflection property.

\begin{proposition} \lll{refell} A beam of light fired from one vertex which
reflects off the ellipse will strike the other vertex. In other
words, $\ang F_1OC = \ang F_2OD$.
\end{proposition}

\includegraphics[width=100mm, height=70mm]{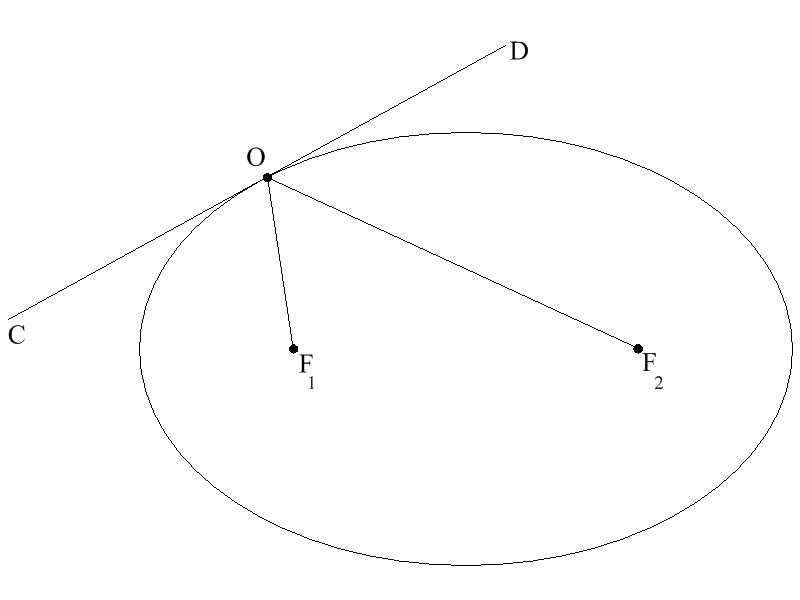}

{\bf Proof:} Draw the circle with center $O$ which is tangent to the
two original circles. If point $O$ moves infinitesimally along the
ellipse, the circle with $O$ at the center will expand or contract.
Since we are only moving infinitesimally, we can replace the two
original circles with their tangents at points $A$ and $B$ and apply
Lemma \ref{tangent}. These tangents are also tangent to the circle
with $O$ at the center, and are therefore perpendicular to
$\ova{AO}$ and $\ova{BO}$. Applying Lemma \ref{tangent}, we conclude
that a tangent vector to the ellipse at $O$ is given by $\ova{AO} +
\ova{BO}$, and thus the tangent bisects $\ang BOA$ so that $\ang BOC
= \ang AOC = \ang F_2OD$ \qed

\includegraphics[width=100mm, height=70mm]{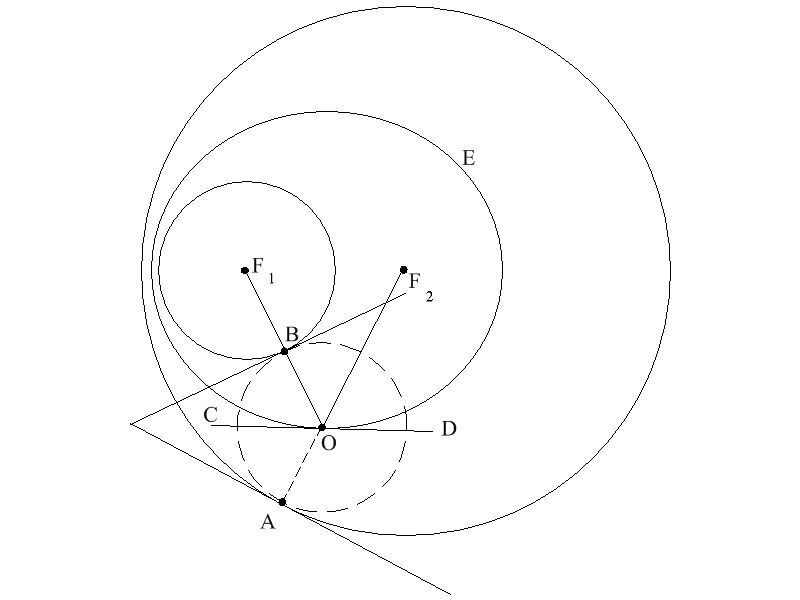}

{\bf Remarks:} i)If the two foci coincide, then the curve is a
circle, and we see that circles are just special cases of ellipses.

ii) The property of ellipses given in Proposition \ref{sum} is often
given as the defining property of ellipses. It is left to the reader
to show that any curve that satisfies Proposition \ref{sum} can be
created by the construction that we have used to define ellipses.

iii) All ellipses besides circles can be created by our construction
where the smaller circle is internally tangent to the larger, and it
is generally simpler to assume that the two circles are tangent. We
didn't want to to assume that here, though, as we want to include
circles. For the remaining conics we lose no generality in assuming
tangency, and we will do so, though it is not necessary.

\vski

Suppose now we have two circles of different radii which are
externally tangent. Let $H$ be the set of all points which are
centers of circles which are externally tangent to both circles. $H$
is a curve, and this curve is known as a {\it hyperbola}. The
centers of the original circles, $F_1$ and $F_2$, are known as the
foci of the hyperbola

\includegraphics[width=100mm, height=70mm]{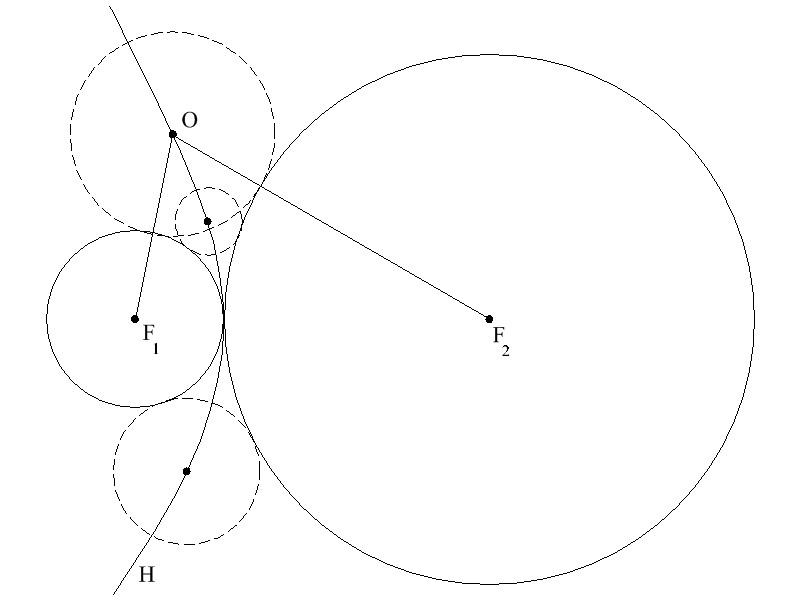}

\begin{proposition} The difference between the distances from the foci to any point on the hyperbola
is a constant.
\end{proposition}

{\bf Proof:} Let $r_1$, $r_2$, and $r_o$ be the radii of the circles
centered at $F_1,F_2$, and $O$ in the above picture. Then

\be F_2O-F_1O = (r_2+r_o) - (r_1 + r_o) = r_2-r_1 \ee
This is a constant, so we are done. \qed

As with the ellipse, we have a very pretty reflection property.

\begin{proposition} Suppose we fire a beam of light from infinity(in other
words, from outside the picture) at one of the foci. If the beam
strikes the hyperbola before reaching the focus, it will reflect off
the hyperbola and strike the other focus. That is, $\ang F_2OY =
\ang \infty OX$.
\end{proposition}

\includegraphics[width=100mm, height=70mm]{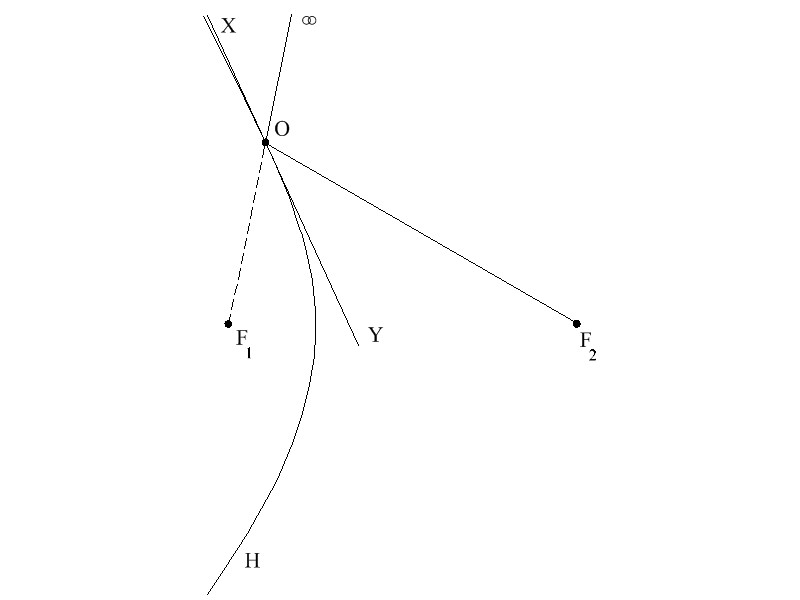}

{\bf Proof:} Pick a point $O$ on the hyperbola, and draw all
relevant circles. We must show $\ang F_1OY = \ang F_2OY$. But this
follows directly from Lemma \ref{tangent} as in Proposition
\ref{refell}. \qed

\includegraphics[width=100mm, height=70mm]{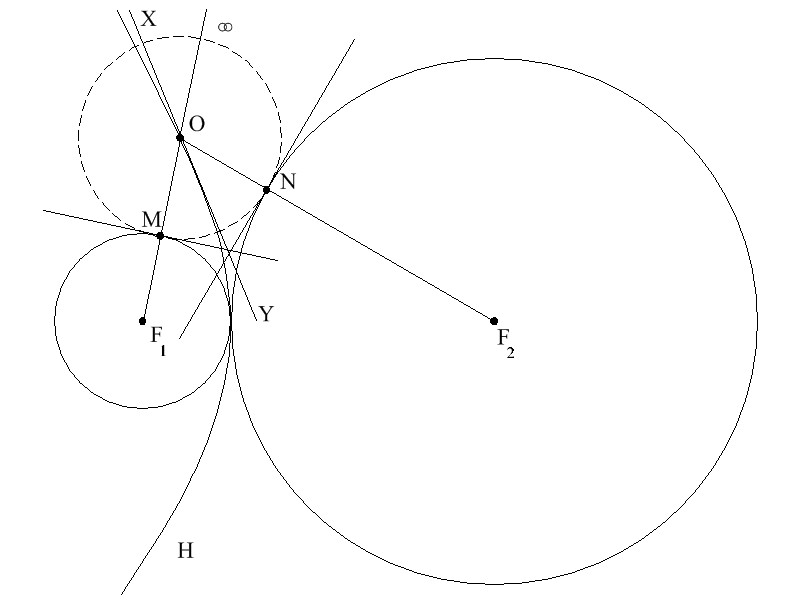}

{\bf Remarks:} i) If the two original circles have the same radii,
the hyperbola reduces to a straight line.

ii) In most books on conic sections, the hyperbola consists of two
parts, one such as we have described, the other the mirror image
reflected around the midpoint of $F_1F_2$. This mirror image will be
obtained by the same construction, interchanging the radii of the
circles centered at $F_1$ and $F_2$.

\vski

The last conic section is the parabola. Suppose that we have a
circle of radius $R$ tangent to a straight line. Let $P$ be the set
of all points which are the centers of circles tangent both to the
line and externally to the circle. $P$ is a curve, known as the
$parabola$. The center of the original circle is called the {\it
focus}.

\includegraphics[width=100mm, height=70mm]{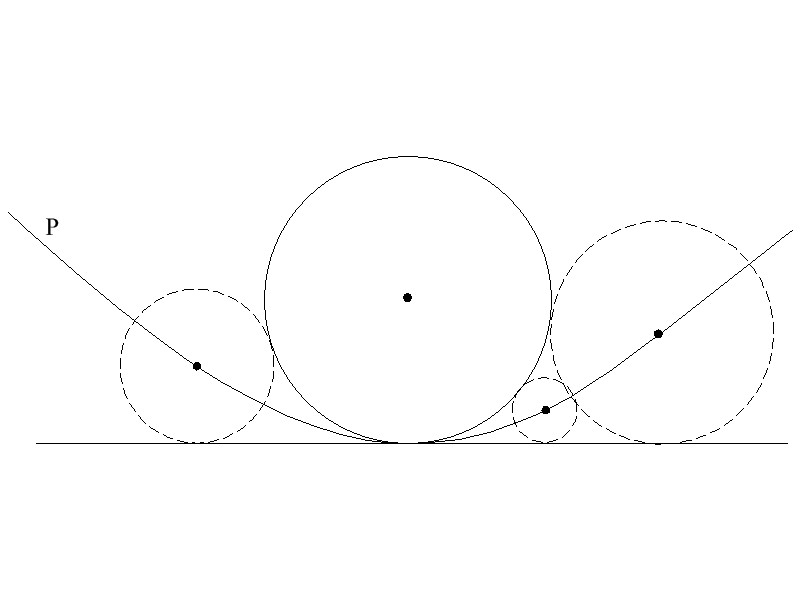}

Displacing the original line $R$ units down forms a new line $L$,
known as the {\it directrix}.

\begin{proposition} The distance from any point on the parabola to the
focus is equal to the distance from that point to the directrix.
\end{proposition}

\includegraphics[width=100mm, height=70mm]{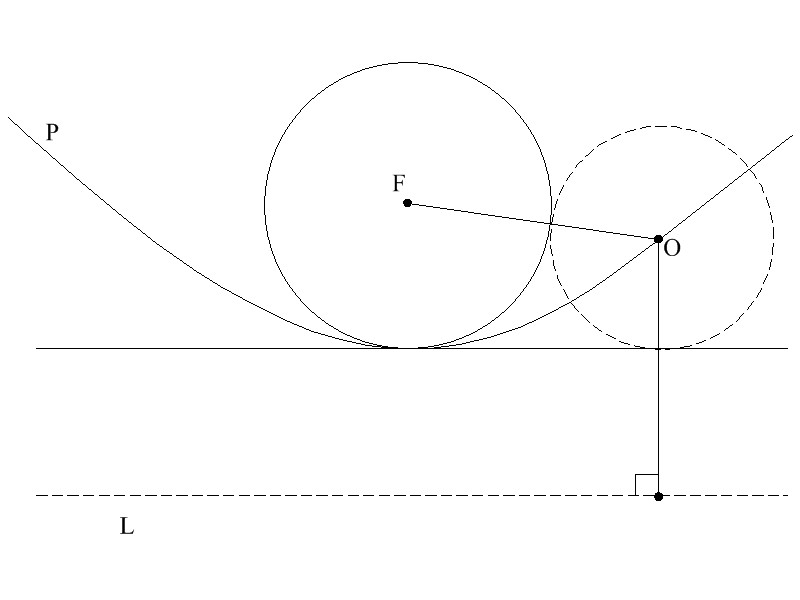}

{\bf Proof:} Let $O$ be a point on the parabola, and let the
corresponding circle have radius $R_O$. Then the distance from the
focus to $O$ is $R + R_o$, as is the distance from the directrix to
$O$. \qed

Finally we have the reflection property, the one that Archimedes is
reputed to have used to torch attacking Roman ships.

\begin{proposition} A beam of light coming straight down will strike
the parabola, reflect off, and hit the focus. In other words, in the
picture below $\ang qpy = \ang opx$.
\end{proposition}

{\bf Proof:} Using Lemma \ref{tangent} as in Proposition
\ref{refell}, $\ang xOa = \ang xOb = \ang qOy$. \qed

\includegraphics[width=100mm, height=80mm]{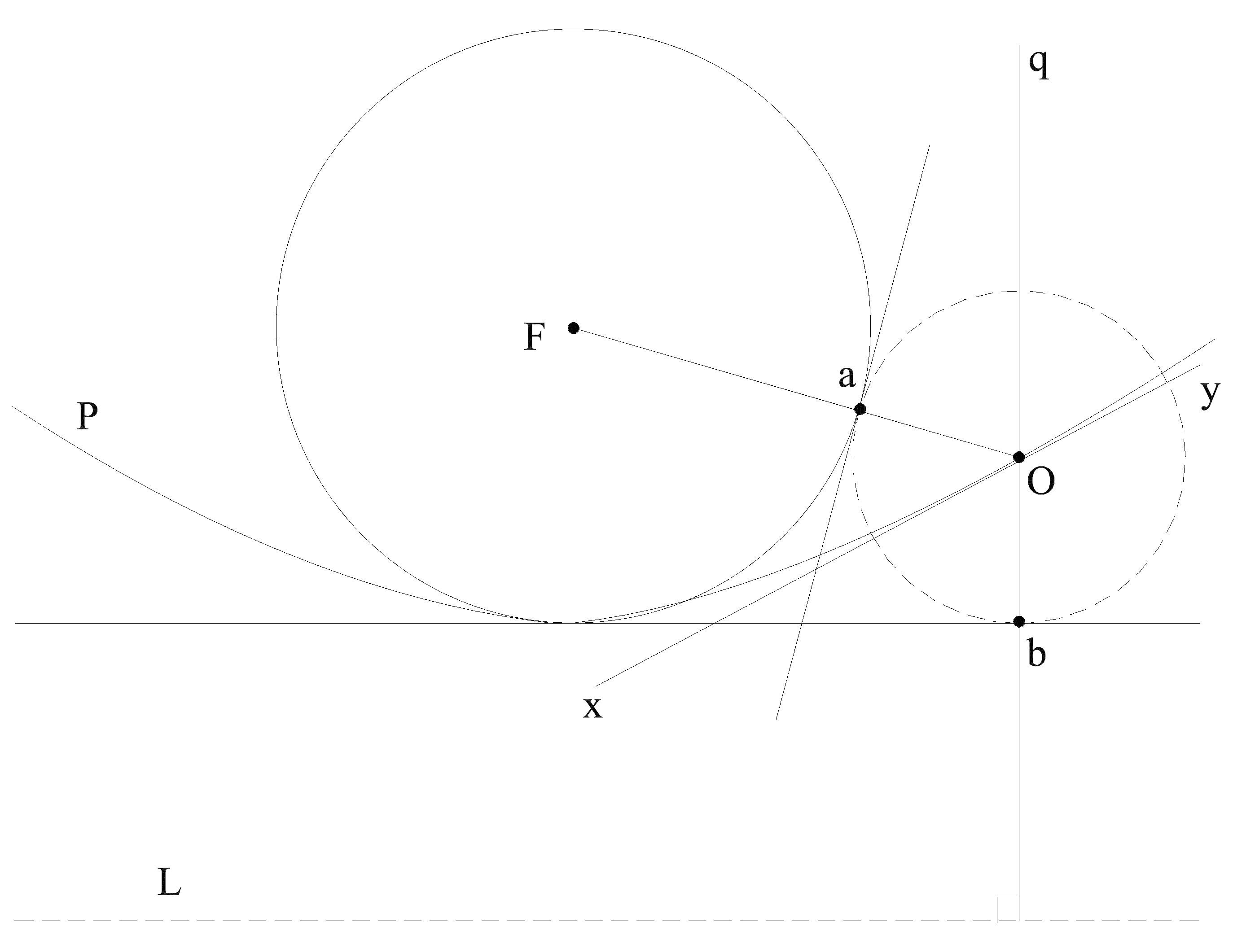}

A plane which cuts through a
cone without touching the vertex gives one of the conic sections.
The slope of the plane determines which conic we obtain. A plane
whose slope is less than the slope of the cone produces an ellipse.

\includegraphics[width=100mm, height=80mm]{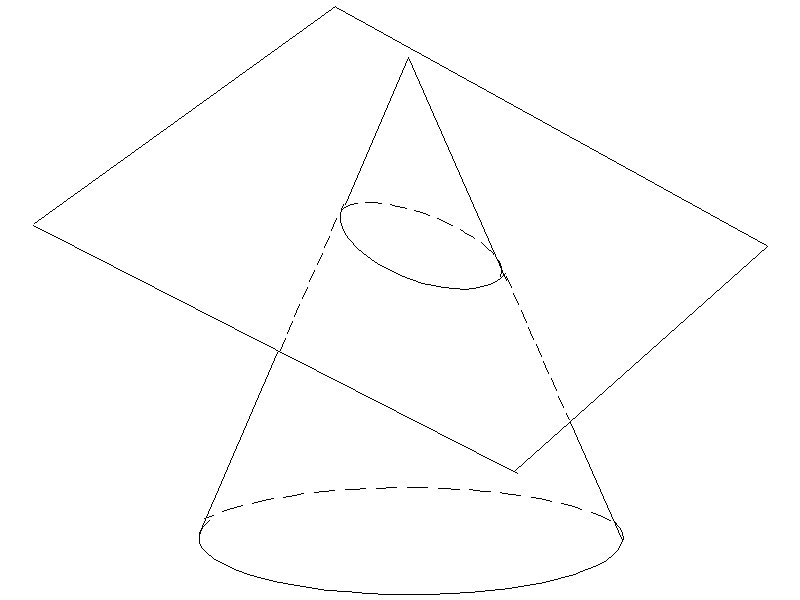}

A plane with a greater slope than the slope of the cone gives a
hyperbola.

\includegraphics[width=100mm, height=80mm]{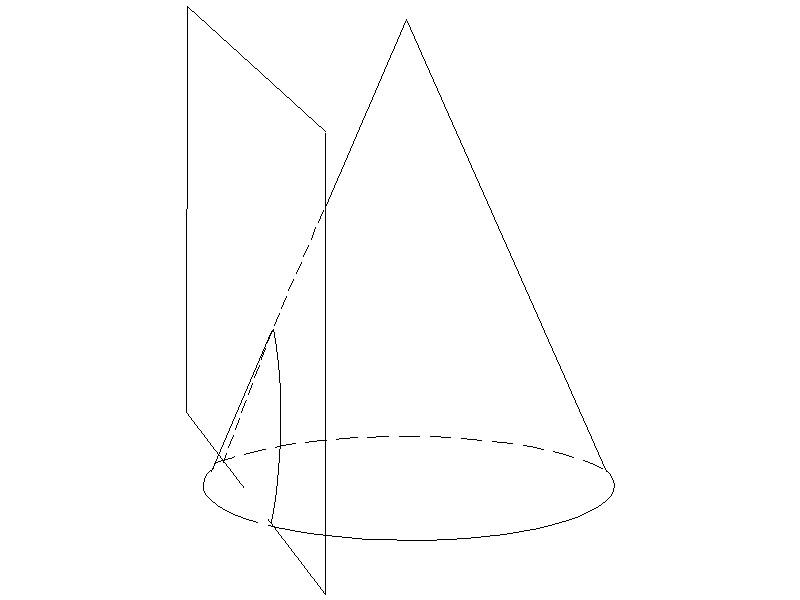}

A plane whose slope is the same as the slope of the cone gives a
parabola.

\includegraphics[width=100mm, height=80mm]{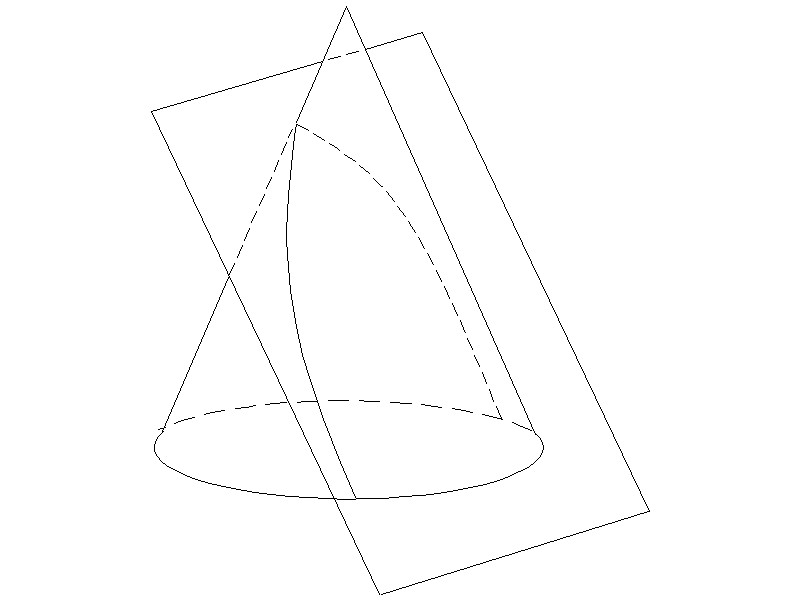}

Let us examine why this is so, beginning with the
ellipse. Suppose that we have a sphere $S$ and a point $O$ away from
the sphere in three dimensions. Then the set of all rays beginning
at the point which are tangent to the sphere will form a cone, and
we will say that the sphere is {\it inscribed} in the cone. The set
of points on the surface of $S$ which are touched by the cone will
form {\it the circle of
tangency}. The points on this circle are equally distant from $O$,
since rotating the entire picture around a vertical axis through $O$
by any angle does not change the picture.

\includegraphics[width=100mm, height=80mm]{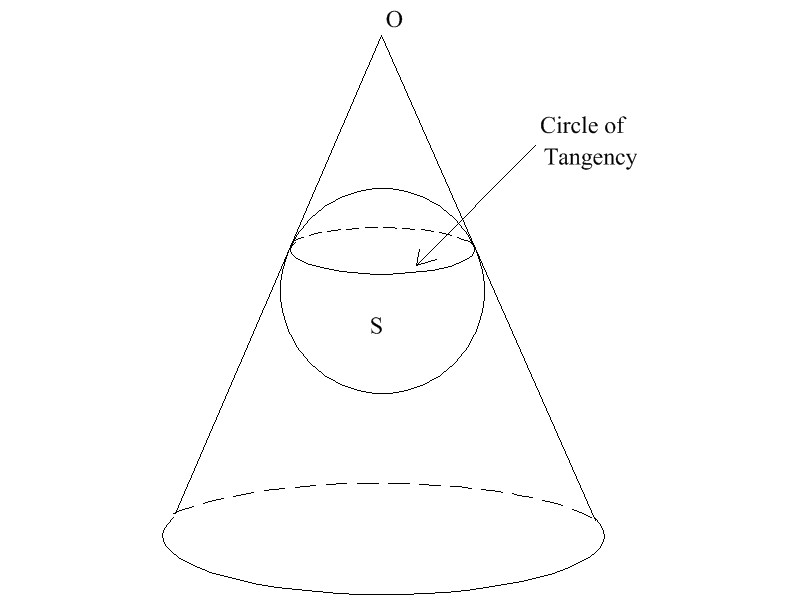}

Suppose we slice a cone with a plane $P$ as shown below. Let
a small sphere be inscribed in the cone above the ellipse, and a
large sphere inscribed below the ellipse. Expand the small sphere
while keeping it inscribed in the cone until it is tangent to $P$ at
point $a$. Contract the large sphere in the same way until it is
tangent to $P$ at $b$. The result is that $a$ and $b$ are the foci
of the ellipse. To see this, let $q$ be a point on the ellipse. Let
$C_1$ and $C_2$ be the circles of tangency of the two spheres, and
let $p_1$ and $p_2$ be the points on $C_1$ and $C_2$ closest to $q$,
i.e. such that $p_2qp_1O$ is a straight line. Since $qb$ and $qp_2$
are both tangent to the larger sphere, their lengths are equal. The
same is true for $qa$ and $qp_1$. Thus, $qb + qa = p_1p_2$, which is
a constant independent of the choice of $q$. We see that the curve
formed is indeed an ellipse, with $a$ and $b$ as its foci.

\includegraphics[width=100mm, height=80mm]{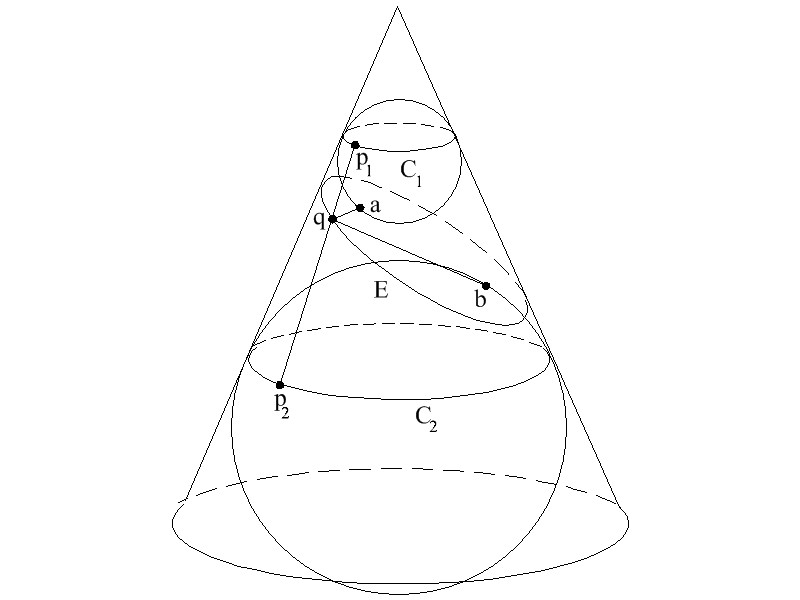}

Now let's do something similar with the hyperbola. We need to extend
the rays through $O$ to form another cone above. Intersect a plane
$P$ with the cone to form a curve $H$, as shown below. Start with
two small inscribed spheres above and below $O$, and expand them
until they are tangent to $P$ at points $a$ and $b$. Let $C_1$ and
$C_2$ be the circles of tangency of the spheres. Choose a point $q$
on $H$, and let $p_1$ and $p_2$ be the points on $C_1$ and $C_2$ so
that $q p_2 O p_1$ is a straight line. Since $q b = q p_2$ and $q a
= q p_1$, we see that $q a - q b = p_1 p_2$, which is a constant.
Thus, $H$ is a hyperbola, with $a$ and $b$ the foci.

\includegraphics[width=120mm, height=80mm]{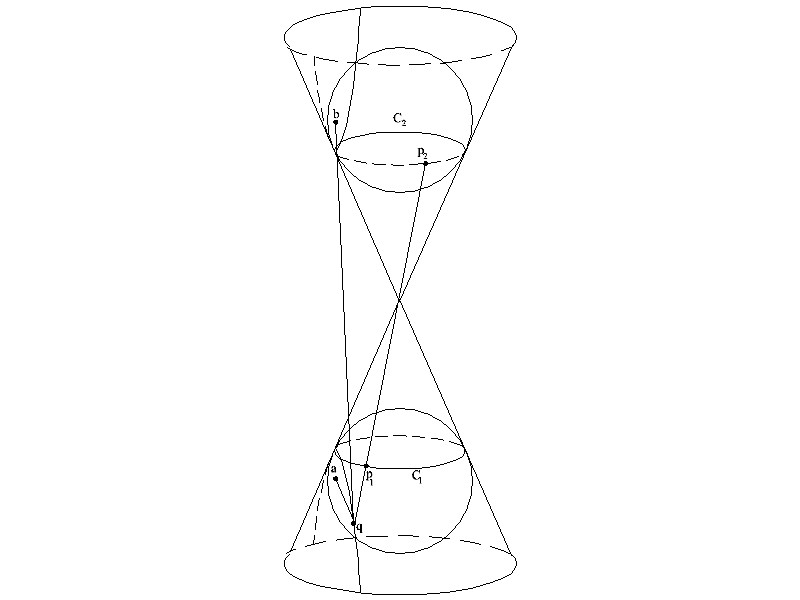}

Now for the parabola. We let $P$ be a plane with the same slope as
the cone, as shown below. Expand a small inscribed sphere near $O$
until it becomes tangent to $P$ at $a$. Let $C$ be the circle of
tangency of this sphere, and draw a line $L$ on $P$ at the same
height as $C$. If $q$ is a point on the curve, $q a = q p$, where
$p$ is the point on $C$ such that $q pO$ is a straight line. Drop a
perpendicular to line $L$ from $q$ to point $l$. Since $P$ is at the
same angle to the vertical as the side of the cone and $C$ is at the
same height of $L$, $q p = q l$, so that $q a = q p$. We see that
the curve is a parabola with vertex $a$ and directrix $L$.

\includegraphics[width=100mm, height=80mm]{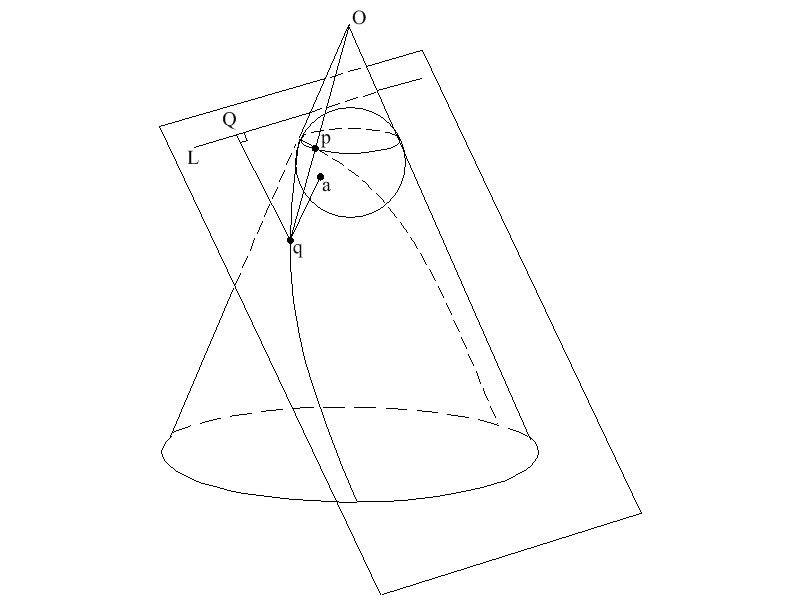}

Before we move on, let us notice one more thing about the parabola example.
Let us keep everything as it was, but change the angle of the plane
$P$. We now have an ellipse or a hyperbola, and it is no longer true
that $q a = q l$. However, it is true that $\frac{q a}{q l} =
\frac{q p}{q l}$ is a constant, since it can be expressed as the
ratio of the slope of the plane and the slope of the cone. We obtain
a new description of the conics.

\begin{theorem} \lll{ecc} Let $L$ be a line and $a$ a point in the plane. The conics can be realized as the set of all
points $p$ such that $\frac{pl}{pa} = e$, where $l$ is the point
closest to $p$ on $L$, and $e>0$ is a constant. If $e < 1$, the
conic is an ellipse. If $e = 1$, the conic is a parabola. If $e>1$,
the conic is a hyperbola.
\end{theorem}

\includegraphics[width=100mm, height=80mm]{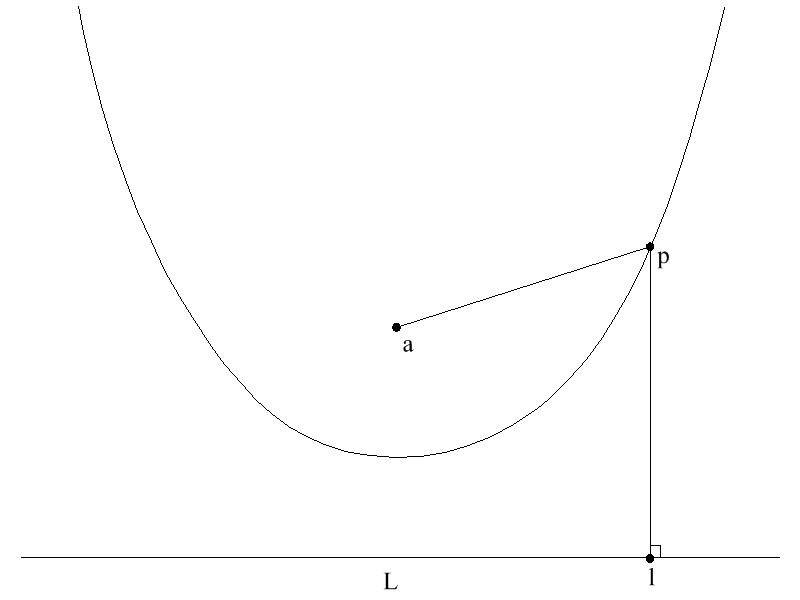}

We refer to $L$ as the {\it directrix}, as has already been
mentioned with the parabola. Let's rip through a few more
propositions on conics.

\begin{proposition} \lll{polar} Let $C$ be a conic with directrix $L$, focus $a$,
and eccentricity $e$. Choose $b$ and $c$ on $C$ such that $b a c$ is
parallel to $L$. Let $p$ be any point on the conic, and let $r$ be
the length of $ap$. Let $\th$ be the angle of $ap$ above $ac$. Then

\be \lll{garb} r = \frac{e(aO)}{1-e\sin \th} \ee

\noi where $O$ is the point on $L$ closest to $a$.
\end{proposition}

\includegraphics[width=100mm, height=80mm]{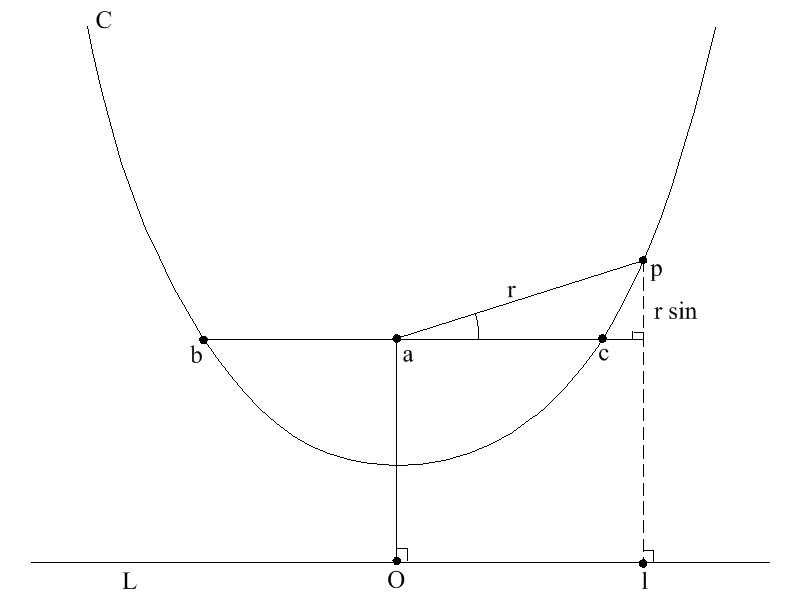}

{\bf Proof:} Drop perpendiculars from $p$ and $a$ to $L$, meeting
$L$ at $l$ and $O$. Then $pl = aO + r\sin \th$. We have

\bea \lll{} \frac{pa}{pl} = e \\ \nn \frac{r}{aO + r \sin \th} = e
\eea
This last equation can be converted to (\ref{garb}) by algebra. \qed

The following is a beautiful little fact which may not be well
known, although it appears in several older books on conics.

\begin{proposition} \lll{right}Let $C$ be a conic with directrix $L$ and focus $a$.
Let $p$ be a point on $C$, and let $T$ be the tangent to $C$ at $p$.
Let $l$ be the point of intersection of $T$ and $L$. Then $\ang pal$
is a right angle
\end{proposition}

\includegraphics[width=100mm, height=80mm]{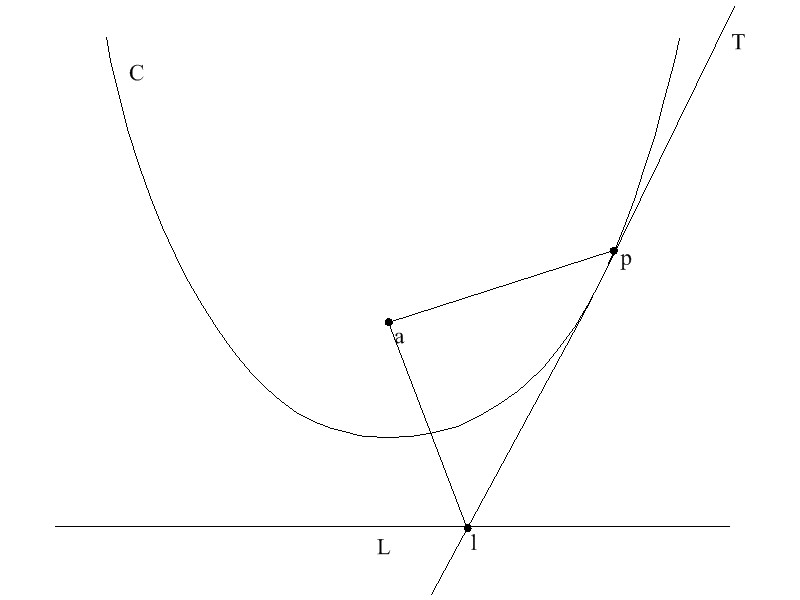}

{\bf Proof:} Let $q$ be a point close to $p$, as shown below. Extend
$pa$ and $q a$ to points $p'$ and $q'$ on $C$, and drop
perpendiculars from $p$ and $q$ to $l_p$ and $l_q$ on $L$. Extend
$pq$ to meet $L$ at $l$, and drop perpendiculars from $l$ to
$\bar{p}$ and $\bar{q}$ on $pp'$ and $qq'$.

\includegraphics[width=100mm, height=80mm]{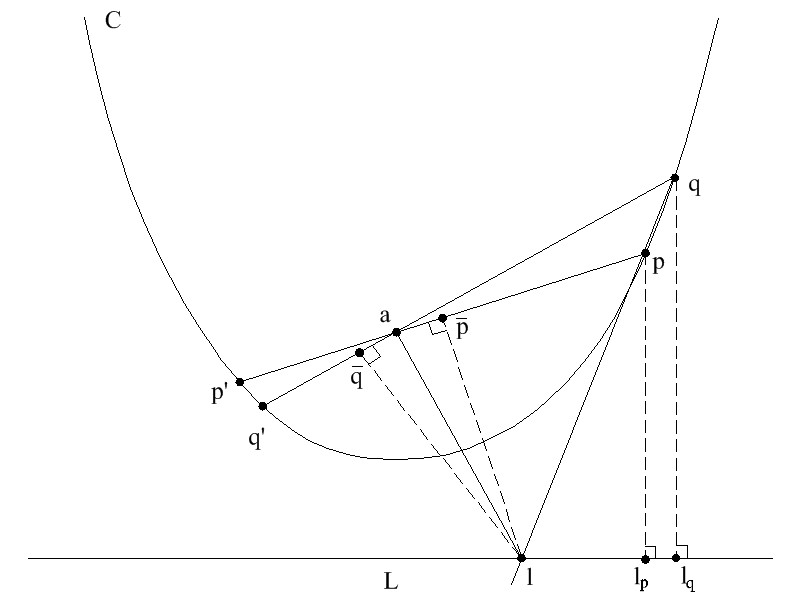}

The first order of business is to show that $\ang q'al = \ang pal$.
We have the following string of equalities.

\be \lll{}\frac{\triangle lap}{\triangle laq} = \frac{lp}{lq} =
\frac{pl_p}{q l_q} = \frac{ap}{aq} \ee
The first equality is because $\triangle lap$ and $\triangle laq$
share a common height to vertex $a$, so $\frac{\triangle
lap}{\triangle laq}$ is equal to the ratio of their bases. The
second is because $\triangle lpl_p \sim \triangle lql_q$. The third
is due to Theorem \ref{ecc}. However, $\triangle lap = (1/2) (ap)(l
\bar{p})$, and $\triangle laq = (1/2) (aq)(l\bar{q})$. We conclude
that $l \bar{p} = l \bar{q}$, and it follows that $\triangle
la\bar{p} \sim \triangle la\bar{q}$. Thus, $\ang q'al = \ang pal$.
Now, hold $p$ fixed and let $q$ approach $p$. Then $q pl$ becomes
the tangent at $p$, and $\ang paq'$ becomes a $180^\circ$ angle. As
this happens, $\ang pal$ becomes a right angle, and we are done.
\qed

The next theorem doesn't have a particularly exciting statement,
but it will be crucial when we begin calculating the orbits of
planets later on.

\begin{theorem} \lll{gook} Let $C$ be a conic, with focus $a$, directrix $L$ and point $O$ chosen on $L$ so that $aO$ and $L$
are perpendicular. Let $p$ be any point on $C$, let $r$ be the
length of $ap$, and let $\aaa$ be the angle made by $ap$ and the
tangent to $C$ at $p$.

\be \lll{}\csc ^2 \aaa = \frac{(e^2 - 1)}{e^2(aO)^2}r^2 +
\frac{2}{e(aO)}r \ee
where $e$ is the eccentricity of $C$.
\end{theorem}

{\bf Remark:} Note that there are two possibilities to choose from
for angle $\aaa$. However, if we label them $\aaa_1, \aaa_2$, we see
$\aaa_1 + \aaa_2 = 180^\circ$, so that $\csc \aaa_1 = \csc \aaa_2$.
In other words, it doesn't matter which one we choose.

\medskip

{\bf Proof:} Choose $\aaa$ as shown below, and let $l$ be the
intersection of the tangent at $p$ with $L$. Then $\ang pal$ is a
right angle, by Theorem \ref{right}. Drop a perpendicular from $a$
to point $O$ on $L$. It is clear from the picture that $aO = r\cos
\th \tan \aaa$

\includegraphics[width=100mm, height=80mm]{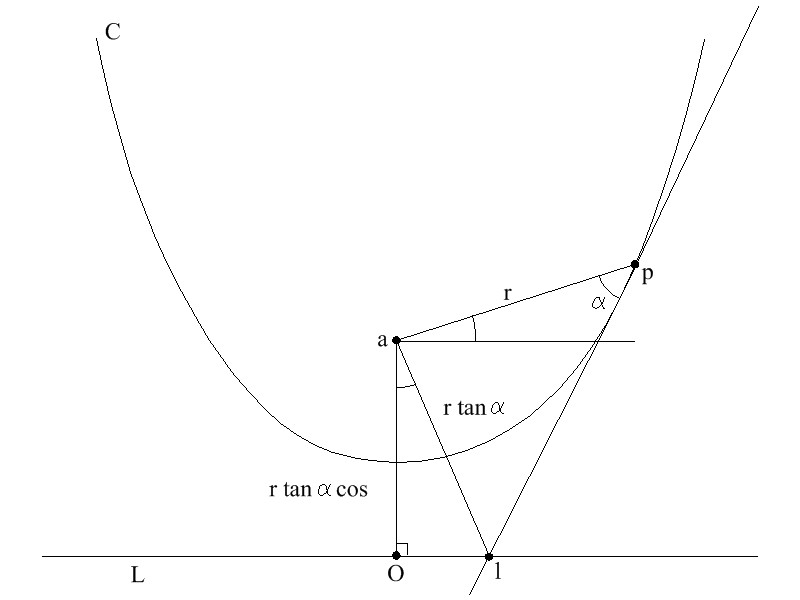}

Squaring this equation gives

\be \lll{ad} (aO)^2 = r^2 \cos^2 \th \tan^2 \aaa \ee

Recall that

\be r = \frac{e(aO)}{1-e\sin \th} \ee

by Theorem \ref{polar}. Rearranging this and squaring gives

\be \sin^2 \th = \frac{(r-e(aO))^2}{r^2e^2} \ee

Thus,

\be \lll{} \cos^2 \th = 1 - \sin^2 \th = \frac{r^2e^2
-(r-e(aO))^2}{r^2e^2} \ee

Plugging this into (\ref{ad}) and rearranging gives

\be \lll{} \tan^2 \aaa = \frac{(aO)^2e^2}{r^2e^2 - (r-e(aO))^2} \ee

Thus,

\be \lll{} \cot^2 \aaa = \frac{r^2e^2 - (r-e(aO))^2}{(aO)^2e^2} =
\frac{r^2(e^2-1)}{e^2(aO)^2} + \frac{2r}{e(aO)} - 1 \ee
Adding $1$ to both sides and using the identity $\csc^2 \aaa = 1 +
\cot^2 \aaa$ completes the proof. \qed

One last proposition about ellipses. {\it Chords} of conics are line segments connecting two points on the conic. The {\it major axis} of an ellipse is the chord passing through the two foci, and the {\it minor axis} is the chord contained in the perpendicular bisector of the major axis.

\begin{proposition} \lll{han} Let $E$ be an ellipse with eccentricity $e$, focus $a$, directrix $L$ and point $O$ chosen on $L$
so that $aO$ and $L$ are perpendicular. Let $Q$ be the center of the
ellipse, and let $X=QB$, $Y=QC$ be the major and minor axes of the
ellipse, respectively. Let $G$ be the area of the ellipse. Then

\be \lll{} e = \frac{\sqrt{X^2-Y^2}}{X} \ee

\be \lll{} (aO) = \frac{Y^2}{\sqrt{X^2-Y^2}} \ee

\be \lll{} X = \frac{(aO)e}{1-e^2} \ee

\be \lll{} Y = \frac{(aO)e}{\sqrt{1-e^2}} \ee

\end{proposition}

\includegraphics[width=100mm, height=80mm]{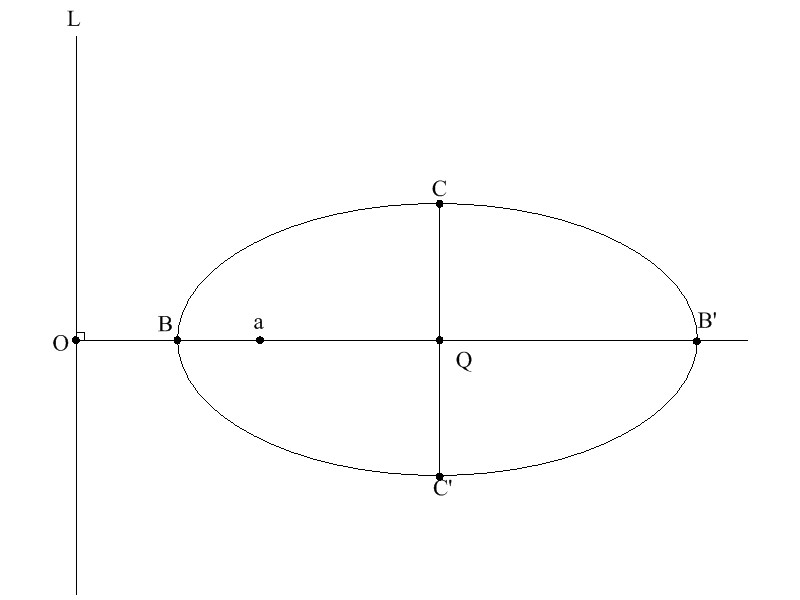}

{\bf Proof:} These relations can be worked out by straightforward
but uninspired calculations in the plane. Happily, there is an
inspired way to do it if we jump to three dimensions. We will
consider an ellipse as the intersection of a cylinder of radius $Z$
with a plane $P$. The same argument as with the cone shows that the
resulting curve is an ellipse; alternatively, we can consider a
cylinder as the limiting case of a cone as we let the height go to
$\infty$ while keeping the base fixed.

\includegraphics[width=100mm, height=80mm]{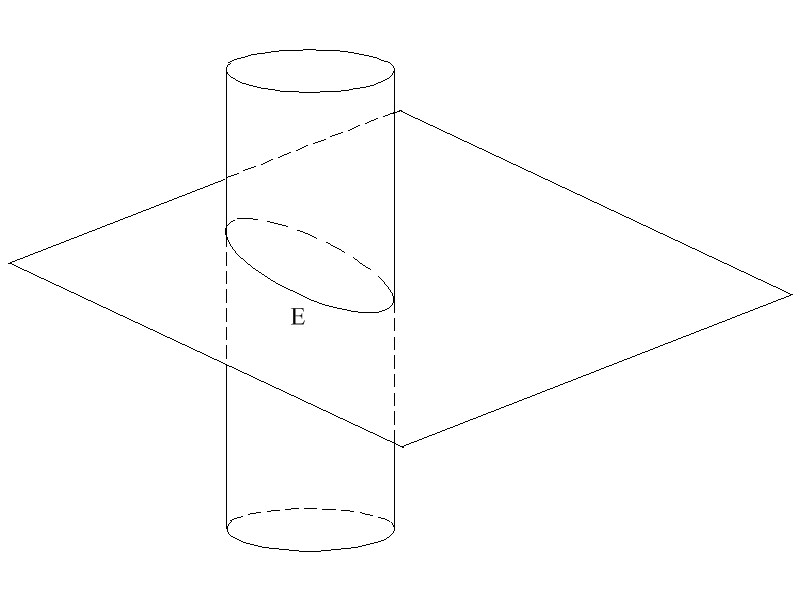}

Inscribe a sphere in the cylinder above $P$, then slide it down
until it is tangent to $P$ at one point. By the same argument as was
used to prove Theorem \ref{ecc}, this point is the focus of the
ellipse, $a$. Furthermore, $L$ is the line on $P$ which is the same
height as the center of the sphere, since the circle of tangency of
the sphere is the circle of that same height. Let us view the entire
setup as a cross section from the side.

\includegraphics[width=100mm, height=80mm]{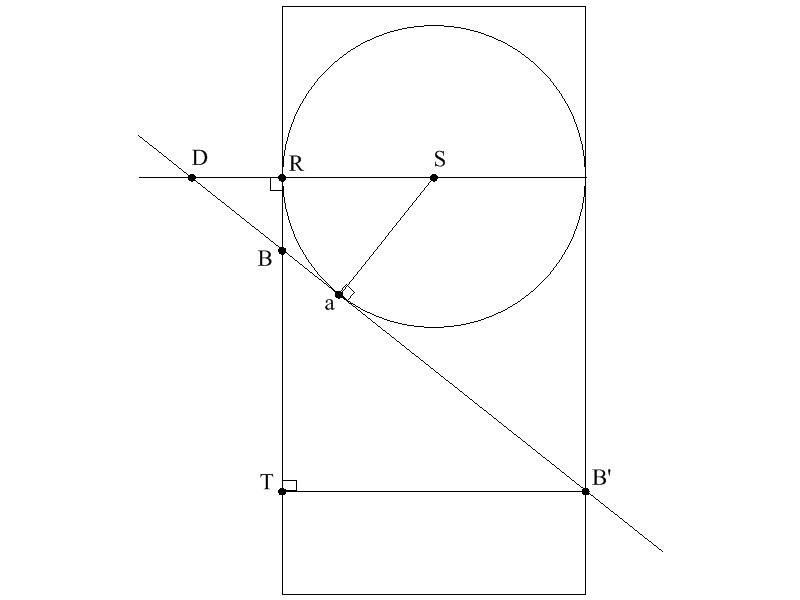}

We have

\be e = \frac{aB}{OB} = \frac{RB}{OB} = \frac{TB}{BB'} =
\frac{\sqrt{X^2 - Y^2}}{X} \ee
The first equality is the definition of $e$, the second is by the
equality of tangents from a point($B$) to a circle, the third is
because triangles $BRO$ and $BTB'$ are similar, and the fourth is
because $BB' = 2X$, $B'T = 2Y$, and therefore $BT = 2\sqrt{X^2
-Y^2}$. Furthermore, triangles $OaS$ and $B'TB$ are similar, so that

\be \lll{} \frac{aO}{Sa} = \frac{B'T}{TB} = \frac{Y}{\sqrt{x^2-Y^2}}
\ee
As $Sa = Y$, we see that $aO = \frac{Y^2}{\sqrt{X^2 - Y^2}}$. We
have established the first two relations. The final two are easy
from this point.

\be \lll{} \frac{(aO)e}{1-e^2} = \frac{(Y^2/X)}{1-(X^2-Y^2)/X^2} = X
\ee

\be \lll{} \frac{(aO)e}{\sqrt{1-e^2}} =
\frac{(Y^2/X)}{\sqrt{1-(X^2-Y^2)/X^2}} = Y \ee

\qed

\section{Kepler, Newton, and experimental data}

Having spent some time in the land of Euclidean geometry, let's return for a bit to the real world. In the early 1600's, Johannes Kepler observed the following rules, which are now famous.

\vspace{.3in}

1. The planets(including the earth) revolve in ellipses about the
sun, with the sun at one of the foci of the ellipse.

\vspace{.3in}

2. A line from the sun to a planet sweeps out equal areas in equal
times.

\vspace{.3in}

3. The time it took for a planet to revolve once around the sun is
proportional to the length of the major axis of the ellipse of
revolution raised to the $3/2$ power.

\vspace{.3in}

Enter Newton. Newton made a few assumptions,
validated through experiments, and proceeded to derive Kepler's laws, or at least come close(there is some debate on whether Newton had something akin to uniqueness for his orbits; see \cite{newtonmag}). We won't need to follow Newton's assumptions to the letter, so we'll make our own. First, let us assume that an object at rest will tend to stay at rest, and any moving object
will tend to continue in a straight line with a constant velocity. A
change in velocity only
happens when the object is acted upon by another object in some way.
Such a change could occur due to a collision, but there are other
ways as well. This property of matter is generally known as {\it inertia}, and was
first proposed by Galileo. To describe our next assumption, suppose that an object $O$ is moving, and that in
the absence of forces acting on it $O$ would move like so in some
amount of time:

\includegraphics[width=100mm, height=80mm]{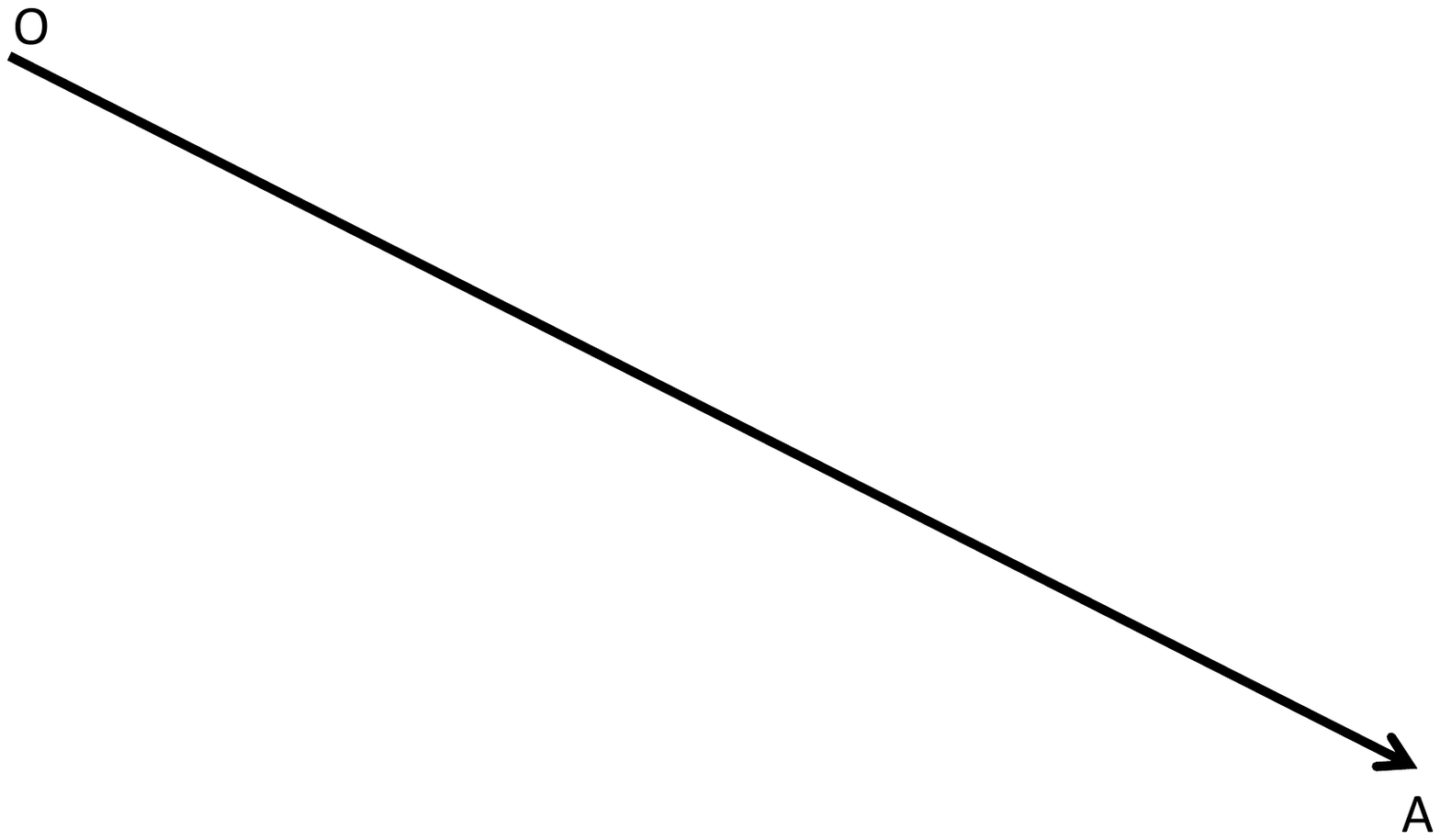}

Suppose in addition that a force acts on $O$, and were this force to act upon $O$
at rest $O$ then in the same amount of time $O$ would move like
so:

\includegraphics[width=100mm, height=80mm]{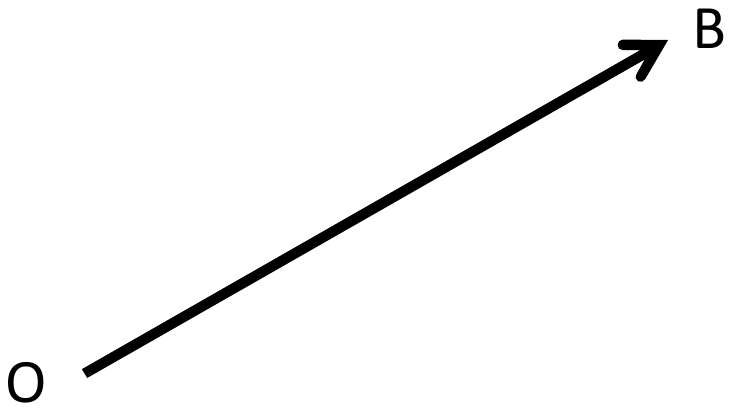}

Then the actual motion of $O$ will be as $Oc$ below:

\includegraphics[width=100mm, height=80mm]{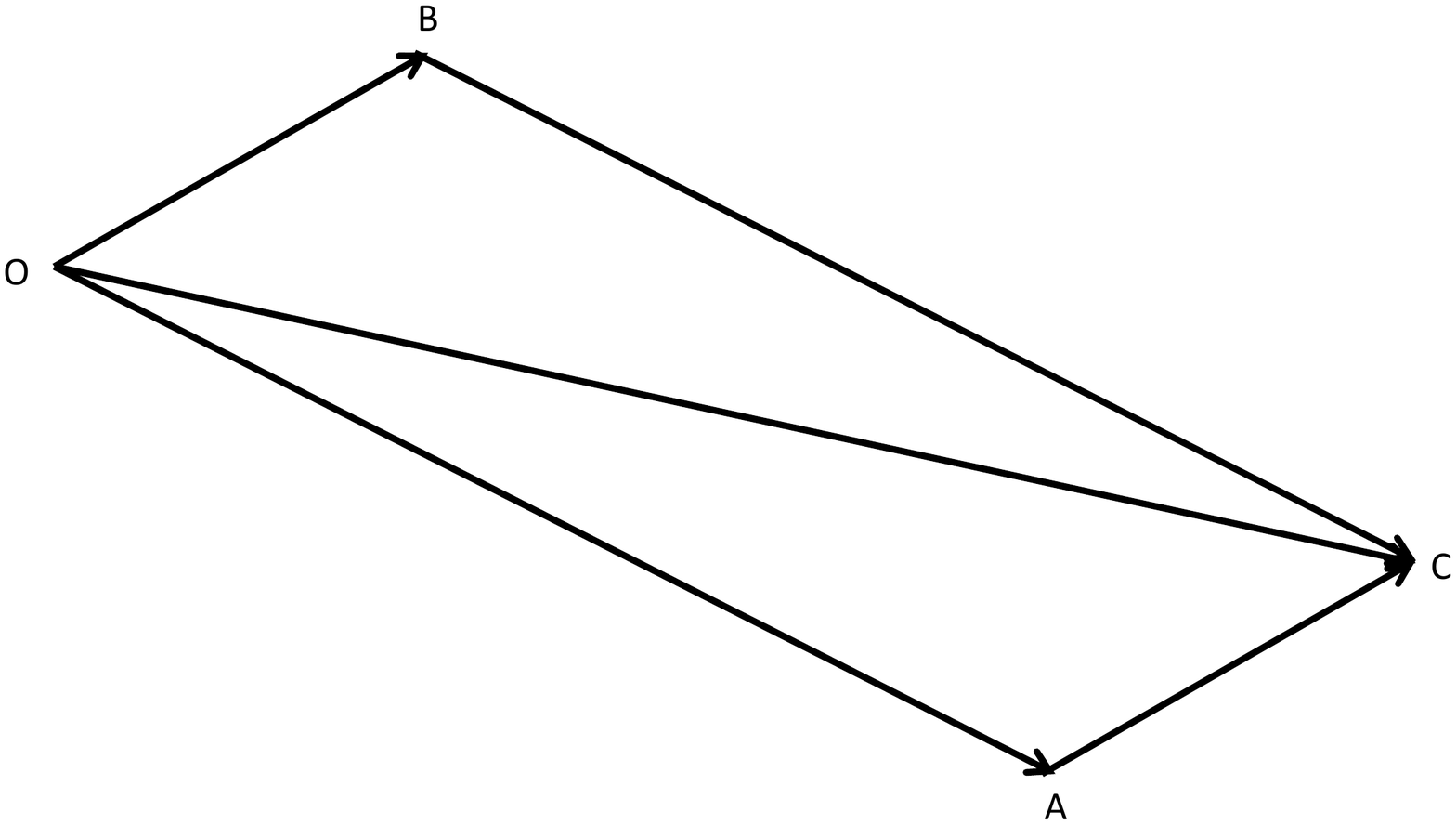}

Furthermore, the same is true if a larger
number of forces acts on an object. In that case we add the
vectors corresponding to each force, and the result gives the movement of
the object. Now, let us see what happens if we suppose that an object moves
according to these assumptions subject only to some force which is
always directed at another object which does not move.

\begin{theorem} \lll{area} Let an object $O$ move through space subject only
to forces directed at a stationary object $S$. Then line segments
connecting $O$ and $S$ will sweep out equal areas in equal times.
\end{theorem}

\includegraphics[width=100mm, height=80mm]{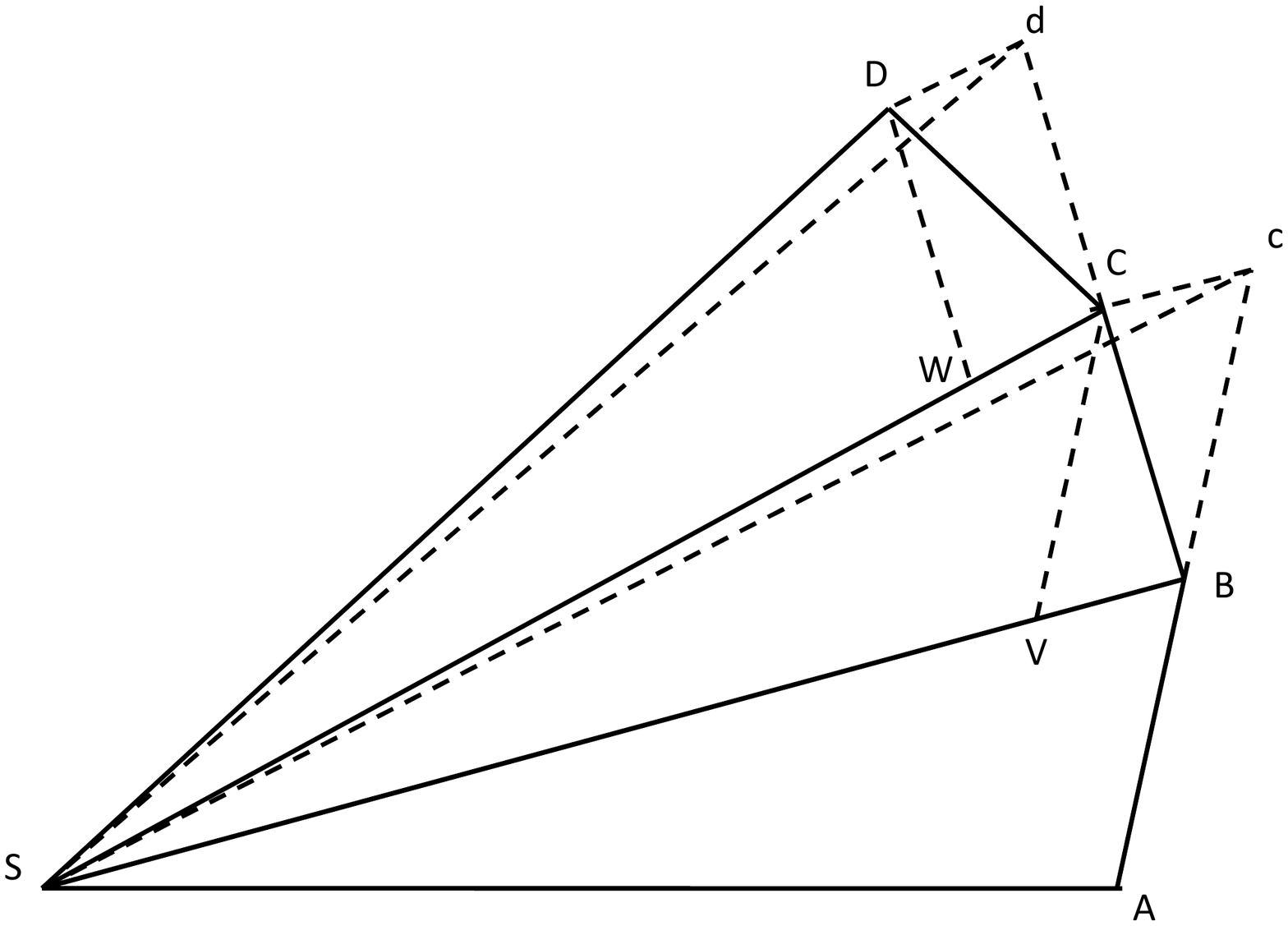}

{\bf Proof:} Let $O$ begin at point $A$, and over some small period
of time $t$ move to point $B$. Then, by the law of inertia, if $O$
were not acted upon by any force it would continue on in a straight
line, arriving at $c$ at time $2t$. However, at point $B$ it is
affected by a "single but great impulse"(Newton's words) directed at
$S$. Let this impulse be equal in magnitude to $BV$. To find the
actual movement of $O$ we complete the parallelogram $BVCc$, and we
see that $O$ resides at point $C$ at time $2t$. Again, in the
absence of force $O$ would move to point $d$ at time $3t$, but we
will assume that $O$ is affected by the impulse $CW$ at point $C$,
so that in fact $O$ resides at $D$ at time $3t$. Now, we must show
that $Area(\ttt SAB) = Area(\ttt SBC) = Area (\ttt SCD)$. Note first
that $Area(\ttt SAB) = Area(\ttt SBc)$, since the bases $AB$ and
$Bc$ of the two triangles are equal, and their common altitude is
the perpendicular dropped from $S$ to the line containing $AB$.
Furthermore, $Area(\ttt SBc) = Area(\ttt SBC)$, since these two
triangles share the base $SB$ and have equal altitudes due to the
fact that $Cc$ is parallel to $SB$. Thus, $Area(\ttt SAB) =
Area(\ttt SBC)$. The same argument shows that $Area(\ttt SBC) =
Area(\ttt SCD)$, and therefore $\ttt SAB, \ttt SBC,$ and $\ttt SCD$
all have the same area.

Now, in general, a force will act continuously and not with isolated
impulses. However, if we let the number of triangles in the above
argument increase to a very large number we will obtain an excellent
approximation to the true path, and in all cases equal areas will be
swept out in equal times. Letting the number of triangles go to
infinity, we obtain a smooth curve in which equal areas are swept
out in equal times, and this gives us the general case. \qed

To steal a line from Richard Feynman(unrelated to this discussion), if a person cannot see the
connection between this theorem and Kepler's second law, then they
have no soul. We have no choice but to at least guess that the
planets move in their orbits due to an acceleration which is always directed at the sun. The
following is a corollary that will be somewhat easier to apply.

\begin{corollary} \lll{prop} Let an object $O$ move through space subject only
to forces directed at a stationary object $S$. Then there is a
constant $k$ such that, for any two points $P$ and $Q$ on the orbit,
we have

\be \lll{sound} T = kA \ee
where $T$ is the time it takes the object to move from point $P$ to
point $Q$, and $A$ is the area of sector $SPQ$.
\end{corollary}

\vspace{0.1in}

{\bf Proof:} Let $A'$ be the total area of the orbit, let $T'$ be
the time it takes the object to complete one orbit, and let $k =
\frac{T'}{A'}$. If we divide the orbit into $N$ pieces of time
$\frac{T'}{N}$, then they must each contain equal areas by Theorem
\ref{area}, and this area must be $\frac{A'}{N}$. Thus, we see that
(\ref{sound}) is satisfied for any time interval which is a rational
multiple of the period of the orbit. An arbitrary interval can be
approximated as closely as we please by a rational interval, so the
result follows. \qed

Let us consider, now, how an object at rest moves under the effect
of a constant continuous force.

\begin{theorem} \lll{accel} Suppose an object $P$ at rest at time $0$ experiences a
constant acceleration $a$. Then the displacement from the initial position at
time $T$ is equal to $(1/2)aT^2$.
\end{theorem}

{\bf Proof:} The velocity at time $t$ is given by $at$. Divide the
interval $[0,T]$ into $N$ equal intervals, $[0,T/N]$, $[T/N,2T/N]$,
... , $[(N-1)T/N, T]$. Let us assume as an approximation that the
velocity over the interval $[(n-1)T/N, nT/N]$ is $anT/N$. Thus, the
distance that $P$ covers over the time interval $[(n-1)T/N, nT/N]$
is $(anT/N)(T/N)$. Adding the distance for each of these intervals,
we get

\be \lll{} \frac{aT^2}{N^2}(1 + 2 + ... + N) = \frac{aT^2N(N+1)}{2N^2} \ee
As $N \lar \infty$, we obtain a better and better approximation, and
$\frac{aT^2}{2N} \lar 0$. We conclude that the distance covered at
time $T$ must be $\frac{aT^2}{2}$. \qed

Now that we have this theorem, let us turn our attention to the
force directed at the sun which we theorize keeps the planets in
their orbits. To begin with, it is a logical guess that the strength
of the acceleration $a$ at any point $P$ in space depends only on the distance from
$P$ to the sun, and not on the direction of $P$ from the sun. We'll
assume that from now on. We will use the notation $A \sim B$ for two
varying quantities $A$ and $B$ to denote the property that
$\frac{A}{B} = M$ for some constant $M$. Since elliptical orbits are too
complicated to deal with without warming up a bit, let's simplify by considering
only circular orbits.

\begin{theorem} \lll{circle} Suppose that an object $P$ can be kept in uniform
circular motion at any radius by an acceleration $a$ directed at an
immovable object $S$ which lies at the center of the orbit. Let $I_R$ be the amount of time $P$ takes
to complete an orbit of radius $R$ around $S$. Suppose that $I_R
\sim R^{3/2}$. Then $a \sim 1/R^2$.
\end{theorem}

{\bf Proof:} An orbit of radius $R$ has circumference $2\pi R$, so
if $v$ is the velocity of the object in this orbit then $I_R =
\frac{2 \pi R}{v}$. Thus, $\frac{2 \pi R}{v} \sim R^{3/2}$, i.e. $v
\sim \frac{1}{\sqrt{R}}$. Let $P$ be at the top of the orbit at a
certain time, and at point $b$ a small amount of time $t$ later.
Drop a perpendicular from $b$ to $c$ on the tangent to the orbit at
$P$. Let $d$ be the point opposite $P$ in the orbit.

\includegraphics[width=100mm, height=80mm]{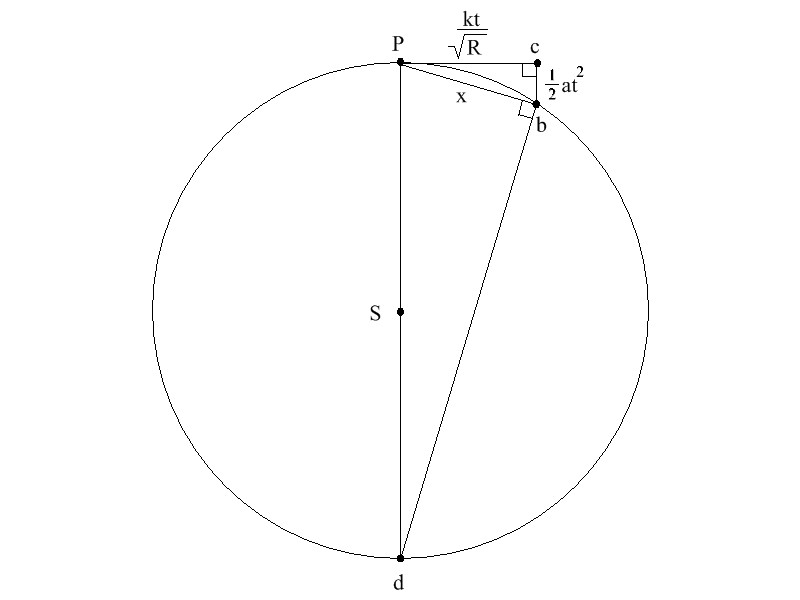}

The movement of $P$ to $b$ can be broken into two components. The
first is $Pc$, which is due to the momentum of the planet at point
$P$ and is therefore approximately equal to $vt$. We know that $v =
\frac{k}{\sqrt R}$, where $k$ is a constant, so $Pc$ is
approximately $\frac{kt}{\sqrt R}$. The other component of the
motion is $cb$, which is due to the acceleration $a$ directed at $S$. We suppose that this acts as a single
impulse at point $P$ directed at $S$, which is why $cb$ is parallel
to $PS$. This impulse has magnitude $\frac{1}{2}at^2$ by the
previous theorem. Since $\ang Pdb$ and $\ang bPc$ subtend equal
arcs, they are equal. Furthermore, $\ang Pcb$ and $\ang Pbd$ are
both right angles. Thus, $\triangle Pdb$ and $\triangle bPc$ are
similar, and

\be \lll{cran} \frac{2R}{x} = \frac{x}{(1/2)at^2} \Rightarrow a =
\frac{x^2}{Rt^2} \ee
Now, we let $b$ approach $P$ so that all of our approximations
become accurate. In doing this $\ang bPc \lar 0$, so that
$\frac{Pb}{Pc} \lar 1$. Thus, we can replace $x$ by $Pc =
\frac{kt}{\sqrt R}$ in (\ref{cran}) to obtain $a = \frac{k^2}{R^2}$.
Thus, $a \sim 1/R^2$. \qed

Eureka! Experimental data in the form of Kepler's third law has led us
to guess that the acceleration keeping the planets in their orbits at any
point is proportional to $1/R^2$, where $R$ is the distance to the
sun. We will say that such an acceleration(or, equivalently, a force) satisfies
an {\it inverse square law}. Now that we have this clue, let's play around with a bit more data to see what happens.
There is an orbit which is closer to us than the orbits of the
planets, namely that of the moon. Let's again approximate this orbit
by a circle. This circle would have a radius of about 385,000 km =
385,000,000 m. One complete orbit takes about 27.3 days $=$
2,358,720 seconds. By the same argument as in the previous theorem,
if $a$ is the acceleration keeping the moon in its orbit, then

\be \lll{} a = \frac{v^2}{R} \ee
$v$ is the circumference over the time, that is $\frac{2 \pi
(385,000,000)}{2,358,720} \approx 1,025 m/s$. Thus,

\be \lll{} a \approx \frac{(1025)^2}{385,000,000} \approx .0027
m/s^2 \ee

So what does this prove? Let's just notice that the force of gravity
at the surface of the earth creates an acceleration of about $9.8
m/s^2$, and that the average radius of the earth is in the ballpark
of $6,367,000 m$. Thus, the ratio between the acceleration at the
surface of the earth and $1/($Distance from surface of earth to
center of earth)$^2$ is about 3.97 x $10^{14}$. Furthermore, the
ratio between the acceleration on the moon and $1/($Distance from
the moon to the center of the earth)$^2$ is approximately 4.00 x
$10^{14}$. Eureka again! These numbers are almost identical, so if
gravity is the acceleration $a$ keeping the moon in orbit, then again we
would have $a \sim 1/R^2$, where $R$ is the distance from an object
to the center of the earth. And if gravity keeps the moon in orbit
around the earth, then why not the planets around the sun? It all
fits together.

There is one troubling objection to this last argument, however,
which the reader may have noticed. That is, I assumed that the force
of gravity somehow originates at the center of the earth, which is a
hard assumption to justify. In fact, in Theorem \ref{circle} this
same assumption was made as well, as we treated the objects as
points without taking into account their size. In the case of Theorem \ref{circle} we can
perhaps argue that the radii of the orbits is so large compared to
the radii of the objects that we may well treat them as points,
but in the more recent argument this is less convincing. How do we
get around this?

Happily, this objection occurred to Newton as well, so we need only
consult {\it Principia}. We begin by assuming that all pieces of matter are
accelerated towards all other pieces of matter in an inverse square law. To be precise,
if the masses
of two small chunks of matter $O_1$ and $O_2$ are $m_1$ and $m_2$, and the distance
between them is $r$, then $O_1$ undergoes an acceleration of $\frac{g m_2}{r^2}$ towards $O_2$, and
$O_2$ undergoes an acceleration of $\frac{g m_2}{r^2}$ towards $O_1$, where $g$ is a (very small) constant.
The ensuing theorem shows that with
spherical objects we can assume that all the mass is concentrated at
the center of the object. Thus, the assumption we made above causes no difficulty.
Before the theorem, let's prove a lemma
that we'll need.

\begin{lemma} \lll{ring} Let $V$ be a very thin ring on the surface of a sphere
centered at a point $a$. Then the surface area of the ring is
approximately $2 \pi x w$, where $x$ is the distance from the ring
to the radius of the sphere through $a$, and $w$ is the width of the
ring.
\end{lemma}

{\bf Proof:} The two pictures below represent the same scenario,
with the second one being a cross section directly from the side.

\includegraphics[width=100mm, height=80mm]{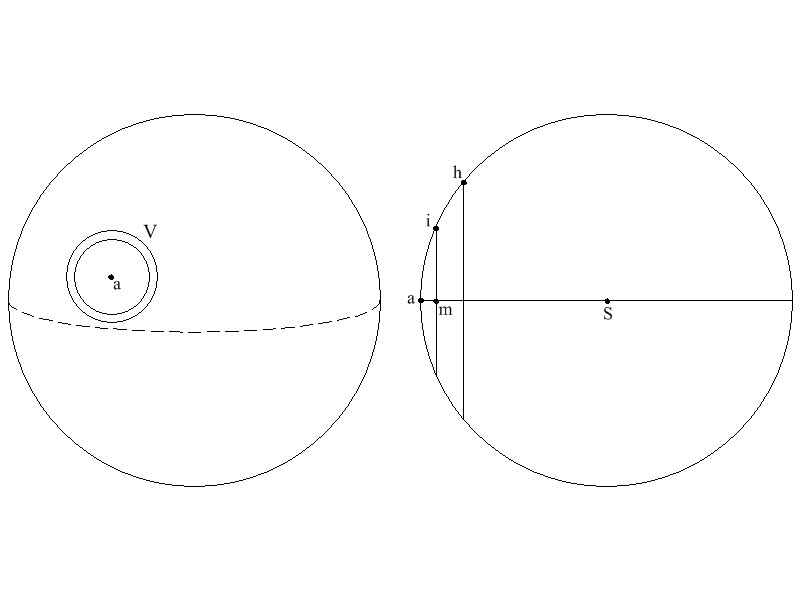}

We want to prove that the area of $V$ is about $2\pi(mi)(ih)$. Since
$ih$ is very small, we can approximate it with a straight line. In
doing so, the ring is approximated by a piece of the top surface of
a cone. Here is the cross section of that cone.

\includegraphics[width=100mm, height=80mm]{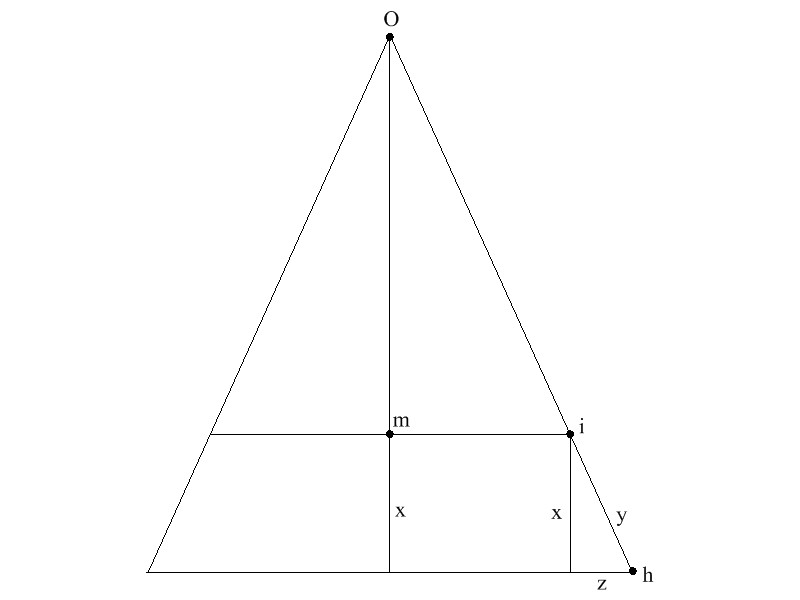}

The area of the top surface of a cone is given by $\pi rs$, where
$r$ is the radius of the base, and $s$ is the distance from the
vertex of the cone to the outside of the base. Thus, in the picture
above, the area of the ring around the cone determined by $ih$ is

\be \lll{} \pi (mi + z)(oi + y) - \pi(mi)(oi) = \pi ((mi)y + (oi)x +
yz) \ee
Now, $y$ and $z$ are both extremely small, so that $yz$ is very
small indeed, much smaller than $(mi)y$ and $(oi)x$. As such, we'll
ignore that term. Furthermore, by similar triangles $\frac{mi}{oi} =
\frac{x}{y}$. We see that the area of the ring is approximately $2
\pi (mi) y = 2 \pi (mi)(ih)$, which is what we set out to prove.
\qed

\begin{theorem}Let $O$ be an unmoving spherical object of uniform mass and center $S$, and let $P$ be a particle
outside $O$. Suppose that $P$ undergoes an acceleration $a$ towards every particle $V$ in $O$,
where $a = \frac{gm_V}{r_V^2}$ with $m_V$ the mass of $V$ and $r_V$ the
distance between $P$ and $V$. Then the overall attraction exerted upon P by
O is inversely proportional to $PS^2$.
\end{theorem}

{\bf Proof:} Let $P$ and $p$ be two identical particles places at
different distances from $O$. Let us suppose first that $O$ is not
actually a solid sphere, but is instead a very thin spherical shell.
We draw two similar diagrams, corresponding to $P$ and $p$, where
the point labels in the first are all capitalized versions of the
lower case labels in the second. These diagrams represent cross
sections of $O$.

\includegraphics[width=150mm, height=120mm]{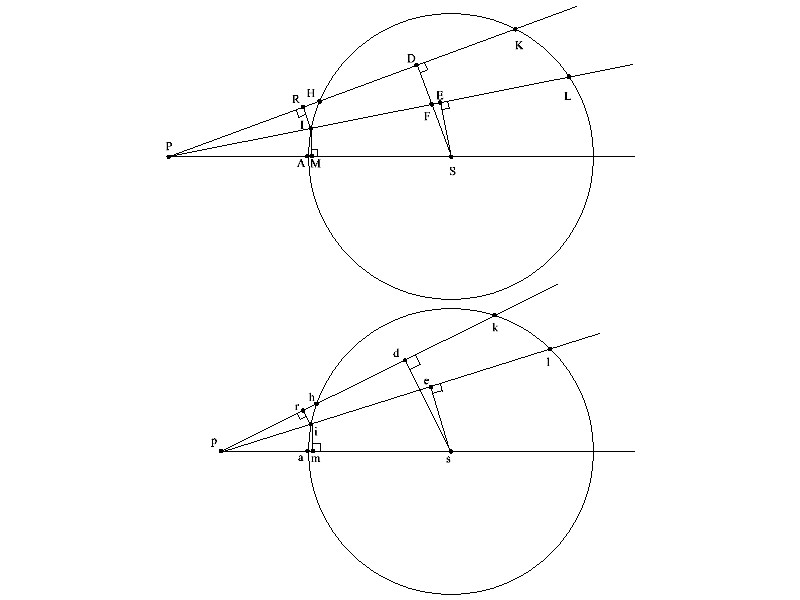}

We will show that the ratio between the accelerations exerted on $P$ and
$p$ by $O$ is proportional to $\frac{ps^2}{PS^2}$. Begin by drawing
a line $PIL$, cutting $O$ at $I$ and $L$, and draw another such
chord $PHK$ such that $\ang IPH$ is very small. On the other
diagram, draw lines $pil$ and $phk$ such that arc $il$ = arc $IL$,
and arc $hk$ = arc $HK$. Let $SE$ and $se$ be perpendicular to $IL$
and $il$, $SD$ and $sd$ perpendicular to $HK$ and $hk$, and $IR$ and
$ir$ perpendicular to $PK$ and $pk$. We will consider the effect on
$P$ due to gravity of the ring created by rotating $HI$ around the
axis $AS$. Given a little piece of this ring, the acceleration on $P$ would
be a constant $k$ times the area of the piece(which is proportional
to the mass of the piece) divided by the distance squared. This
distance between $P$ and the ring will be approximately the length
of $PI$ for any point on the ring, so we'll just use that for the
distance. The total area of the ring is about $(MI)(IH)$ by Lemma
\ref{ring}, but the force exerted by the ring on $P$ is not
$\frac{k(MI)(IH)}{(PI)^2}$. Resolve the acceleration along $PI$ into the
sum of an acceleration along $PM$ and an acceleration along $MI$. The acceleration along
$MI$ will be canceled by a corresponding acceleration from the bottom of
the ring. Thus, all that we need consider is the acceleration along $PM$.
The ratio of this acceleration to the acceleration along $PI$ is $\frac{PM}{PI} =
\frac{PE}{PS}$. Now, since $\ang LPK$ is very small, $PE$ and $PF$
are nearly equal, so we can replace this ratio with $\frac{PF}{PS}$,
and it follows that the total acceleration upon $P$ from the ring formed by
$IH$ is proportional to $\frac{k(MI)(IH)(PF)}{(PS)(PI)^2}$. Of course, the same
argument shows that the acceleration upon $p$ by the ring formed by $ih$ is
proportional to $\frac{k(mi)(ih)(pf)}{(ps)(pi)^2}$. The first thing we need to show
is that

\be \lll{toprove} \frac{ps^2}{PS^2} =
\frac{\Big(\frac{(IH)(IQ)}{(PI)^2}\Big)
\frac{PF}{PS}}{\Big(\frac{(ih)(iq)}{(pi)^2}\Big) \frac{pf}{ps}} \ee
Note that $\frac{PI}{PF} = \frac{RI}{DF}$ by similar triangles, and
similarly $\frac{pf}{pi} = \frac{df}{ri}$. Thus,

\be \lll{} \frac{(PI)(pf)}{(PF)(pi)} = \frac{(RI)(df)}{(DF)(ri)} \ee
We can argue that $DF$ and $df$ are nearly equal, though, as
follows. Arcs $HK$ and $hk$ are equal, as are $IL$ and $il$.
Furthermore, since $\ang LPK$ and $\ang lpk$ are very small, $HK$
and $hk$ are nearly centered in $IL$ and $il$. Thus,

\be \lll{} DF \approx DS - ES = ds - es \approx df \ee
So we can replace $\frac{df}{DF}$ with 1 to get

\be \lll{} \frac{(PI)(pf)}{(PF)(pi)} = \frac{(RI)}{(ri)} \ee
$HI$ and $hi$ are very small, so they are essentially straight
lines. Thus, $\frac{(RI)}{(ri)} = \frac{(HI)\sin (\ang RHI)}{(hi)
\sin (\ang rhi)}$. Now, $\ang RHI$ covers arc $IHK$, and $\ang rhi$
covers arc $ihk$. These arcs are nearly the same, so that $\sin
(\ang RHI) \approx \sin (\ang rhi)$. Thus, we may replace
$\frac{(RI)}{(ri)}$ with $\frac{(HI)}{(hi)}$, and we obtain

\be \lll{rat1} \frac{(PI)(pf)}{(PF)(pi)} = \frac{(HI)}{(hi)} \ee
File this equation away for now. By similar triangles, we have
$\frac{PI}{IM} = \frac{PS}{SE}$, hence $\frac{PI}{PS} =
\frac{IQ}{SE}$. Likewise $\frac{pi}{ps} = \frac{iq}{se}$, but since
$SE = se$ we can write $\frac{pi}{ps} = \frac{iq}{SE}$ instead.
Multiplying the ratios together gives

\be \lll{rat2} \frac{(PI)(ps)}{(PS)(pi)} = \frac{(IQ)}{(iq)} \ee
Now take the product of (\ref{rat1}) and (\ref{rat2}). This gives

\be \lll{}\frac{(PI)^2(pf)(ps)}{(pi)^2(PF)(PS)} =
\frac{(IH)(IQ)}{(ih)(iq)} \ee
Thus

\be \lll{}\frac{(pf)(ps)}{(PF)(PS)} =
\frac{\frac{(IH)(IQ)}{(PS)^2}}{\frac{(ih)(iq)}{(ps)^2}} \ee
And, finally

\be \lll{} \frac{(ps)^2}{(PS)^2} =
\frac{(pf)(ps)(\frac{ps}{pf})}{(PF)(PS)(\frac{PS}{PF})} =
\frac{\Big(\frac{(IH)(IQ)}{(PI)^2}\Big)
\frac{PF}{PS}}{\Big(\frac{(ih)(iq)}{(pi)^2}\Big) \frac{pf}{ps}} \ee
which is what we wanted to show(recall (\ref{toprove})). This proves
that the accelerations on $P$ and $p$ by the rings generated by revolving
$HI$ and $hi$ around $SP$ and $sp$ are in the ratio
$\frac{(ps)^2}{(PS)^2}$. A similar argument shows that the same
holds of the rings generated by $KL$ and $kl$ as well. We can divide
$O$ into a large number of very thin rings. We then have the property
that, for any ring $V$ determined by arc $HI$ there is a unique ring
$v$ determined by the arc $hi$ with $\ang PHI = \ang phi$ and such
that the ratio of the accelerations exerted by $V$ on $P$ and by $v$ on $p$
is $\frac{(ps)^2}{(PS)^2}$. Adding all of the rings together shows
that the ratio of the accelerations exerted by $O$ on $P$ and by $O$ on $p$
is $\frac{(ps)^2}{(PS)^2}$ as well. Recall that this was all done
for a hollow shell $O$. If $O$ is a solid sphere, however, we just
think of it as the sum of a large number of thin, hollow shells. The
ratio for the accelerations from each of the hollow shells on $P$ and $p$
is $\frac{(ps)^2}{(PS)^2}$, and when we add all of them together
that ratio persists. \qed

Having dispensed with that difficulty, let's find others to worry about. The last few theorems give convincing evidence
that there is an acceleration called gravity between all chunks of matter,
inversely proportional to the distance between the chunks squared,
which keeps all of the orbits going. But since I can't just leave
well enough alone, I'm going to pile some more evidence on. Above
we've simplified in every case by assuming circular orbits instead
of elliptical, because ellipses are difficult. Now it's time to take on the elliptical orbit. First,
three pretty lemmas about ellipses.

\begin{lemma} \lll{ell1}Let $E$ be an ellipse with foci $a$ and $b$ and major
axis $Om$. Let $p$ be a point on the ellipse as shown below. Choose
$c'$ on $ap$ and $c$ on $pb$ extended so that $c'Oc$
is parallel to the tangent to the ellipse at $p$. Then $Pc' = Pc =
Om$.
\end{lemma}

\includegraphics[width=100mm, height=80mm]{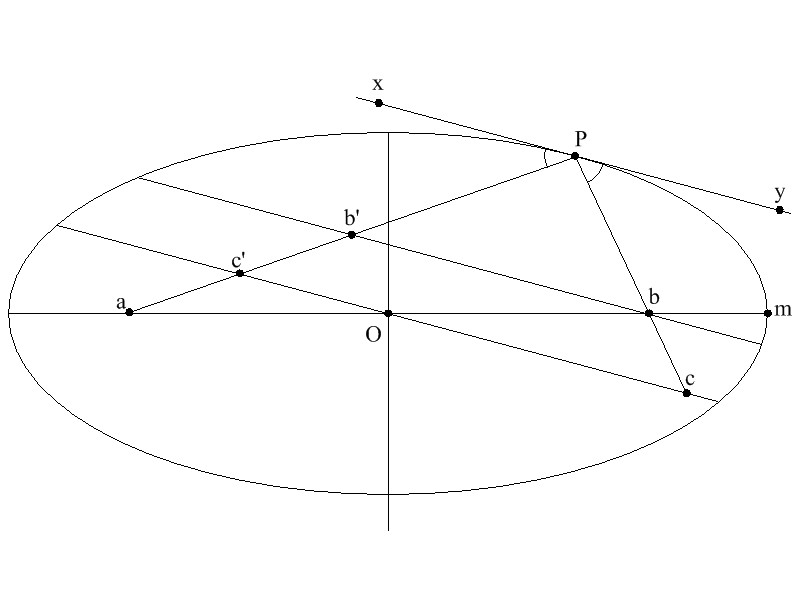}

{\bf Proof:} In the previous section it was shown that $\ang apx = \ang bpy$.
Since $c'Oc$ is parallel to $xpy$, this implies $\ang pc'O = \ang
pcO$, so that triangle $cpc'$ is isoscoles. Thus, $c'p = cp$.
Choose $b'$ on $ap$ so that $bb'$ is parallel to the tangent at $P$.
Triangles $abb'$ and $aOc'$ are similar, so
$\frac{ac'}{c'b'} = \frac{ao}{ob} = 1$. Thus, $ac' = c'b' = cb$. We
have

\bea \nn \lll{}cp + c'p = cb + bp + pc' = c'a + pc' + pb = pa + pb = am + bm = 2
Om
\eea
The second to last equality is due to the fact, proved in the previous section, that
the sum of the distances from the foci of an ellipse to the points on the ellipse is a
constant. The result follows. \qed

\begin{lemma} \lll{ell2} Let $E$ be an ellipse with axes $AC$ and
$BC$. Let $P$ be a point on the ellipse, and let $DK$ be the line
throught the center $C$ of the ellipse parallel to the tangent at
$P$. Let $F$ be on $DK$ so that $PF$ is perpendicular to $DK$. Then
$(PF)(CK) = (AC)(BC)$.
\end{lemma}

\includegraphics[width=100mm, height=80mm]{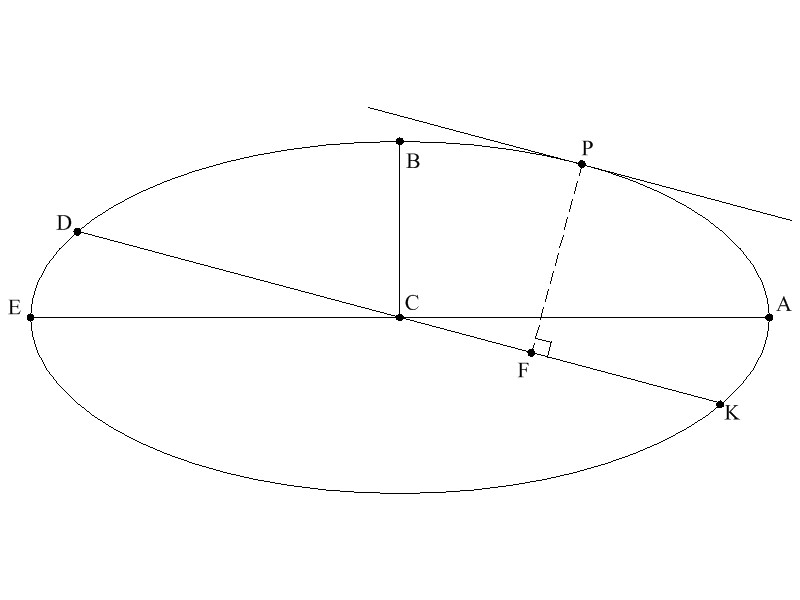}

{\bf Proof:} This is trivial when the ellipse is a circle. For the general case, consider an affine transformation from a circle to the ellipse(if you don't know what an affine transformation is, draw a circle with a marker on a plane of glass and then let sunlight pass through it and strike the ground, casting a shadow of the circle on the ground. By altering the angle of the glass you can produce as the shadow an ellipse of any eccentricity. This is the required transformation). Such transformations preserve ratios of areas, we we can reduce the case of the general ellipse to that of the circle.
\qed

\includegraphics[width=150mm, height=120mm]{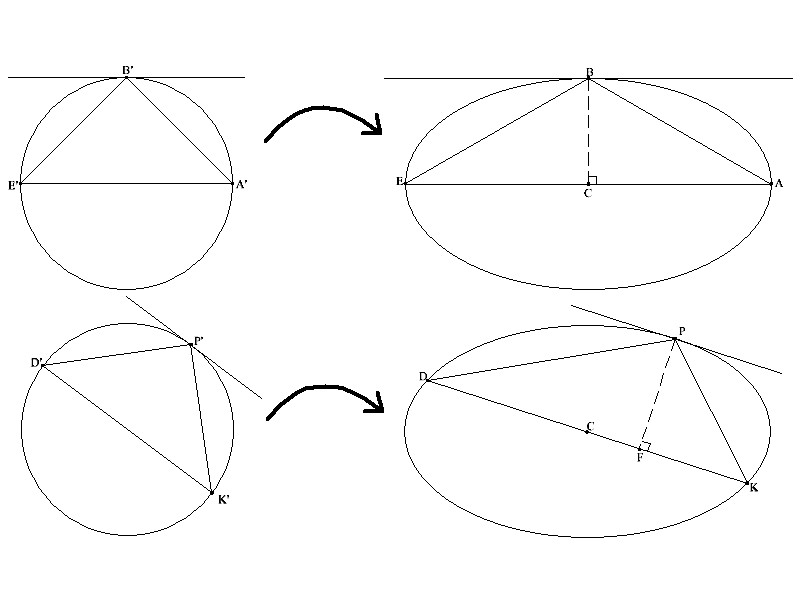}

\begin{lemma} \lll{ell3} Let $GCP$ be a straight line through the center $C$ of
an ellipse. Let $v$ be a point on $GP$ and $Q$ a point on the
ellipse so that $Qv$ is parallel to the tangent at $P$. Draw the
line $DCK$ parallel to the tangent at $P$. Then

\be \lll{ac}\frac{(Gv)(vP)}{(Qv)^2} = \frac{(PC)^2}{(DC)^2} \ee
\end{lemma}
{\bf Proof:} Let us suppose first that the ellipse is in fact a
circle, and let's just assume that $GP$ is vertical for simplicity.

\includegraphics[width=100mm, height=80mm]{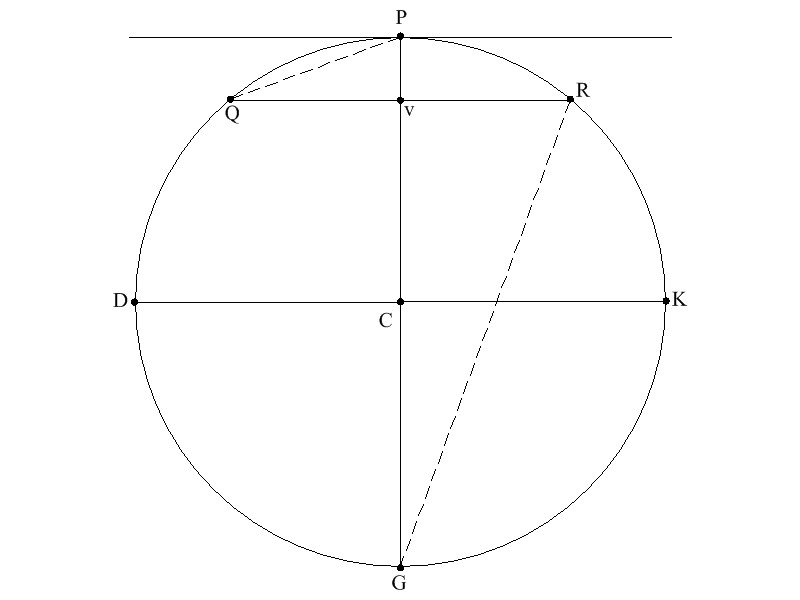}

Clearly $\ang GvR = \ang QvP$, and we also have $\ang GRv = \ang
QPv$ and $\ang vGR = \ang vQP$ as well, since these pairs of angles
cover the same intervals on the circle. Thus, triangles $Cvr$ and
$QvP$ are similar, and $\frac{Qv}{vP} = \frac{Gv}{vR}$, hence
$(Qv)(vR) = (Gv)(vP)$. Since $Qv = vR$, we see

\be \lll{} \frac{(Gv)(vP)(DC)^2}{(Qv)^2(PC)^2} =
\frac{(Gv)(vP)}{(Qv)^2} = 1 \ee
so that (\ref{ac}) holds.

How do we go from a circle to an ellipse? Lesson learned from the
previous lemma, we project the circle with an affine transformation. Choose a circle which projects to the ellipse in
question, and choose points $G',D',Q',P',K',C',v'$ on and in the
circle which project to points $G,D,Q,P,K,C,v$ on and in the
ellipse.

\includegraphics[width=100mm, height=80mm]{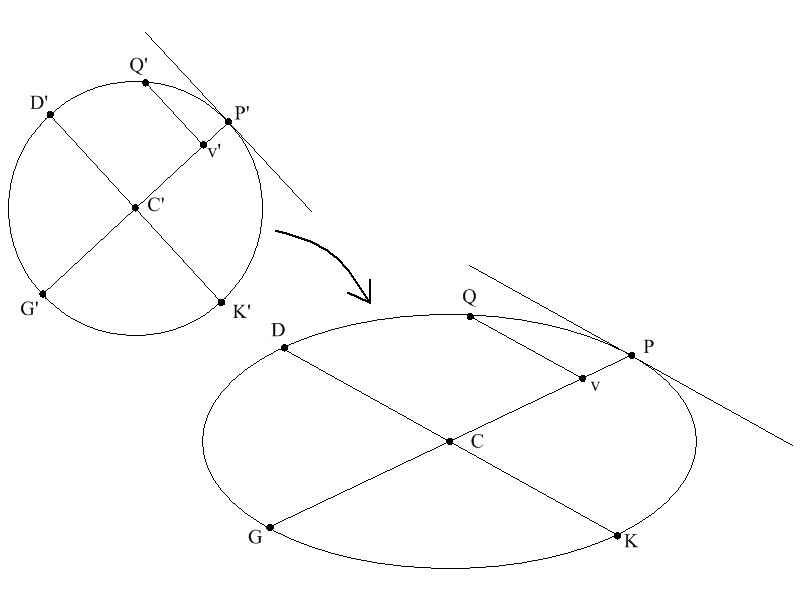}

From the work above, we know that

\be \lll{} \frac{(G'v')(v'P')(D'C')^2}{(Q'v')^2(P'C')^2}  = 1 \ee
A property of this type of transformation is that, if $a'b'$ and $c'd'$ are two line segments which lie on parallel lines and
which project to line segments $ab$ and $cd$, then $\frac{a'b'}{c'd'}=\frac{ab}{cd}$. Thus

\be \lll{} 1 =
\frac{\frac{G'v'}{P'C'}\frac{v'P'}{P'C'}}{\frac{(Q'v')^2}{(D'C')^2}}
= \frac{\frac{Gv}{PC}\frac{vP}{PC}}{\frac{(Qv)^2}{(DC)^2}} \ee
Thus,

\be \lll{ac}\frac{(Gv)(vP)}{(Qv)^2} = \frac{(PC)^2}{(DC)^2} \ee
and we are done. \qed

Now let's look at one of Newton's theorems on elliptical orbits.

\begin{theorem}Suppose that an object $O$ moves in an elliptical orbit due to an
acceleration $a$ towards an unmoving object $S$ which depends only on the
distance $R$ between $S$ and $O$. Then $a \sim \frac{1}{R^2}$.
\end{theorem}

{\bf Proof:} Here is the diagram that appears with this theorem in
{\it Principia}.

\includegraphics[width=100mm, height=80mm]{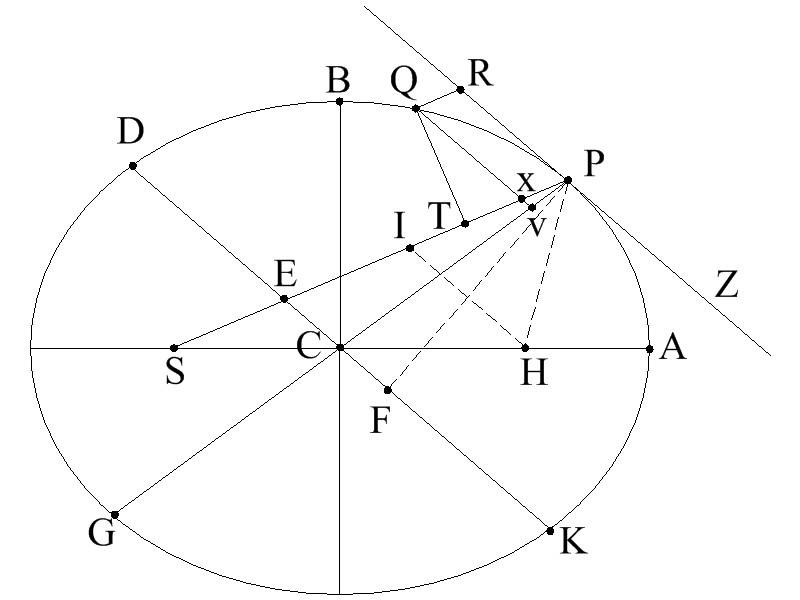}

$S$ and $H$ are the foci of the ellipse. $PF$ is perpendicular to
$DK$, and $QT$ is perpendicular to $PS$. $O$ resides at $P$, and $Q$
is a point very close to $P$. $PF$ is perpendicular to $DK$. In the
style of Newton, we are going to collect a bunch of relations
between ratios, then multiply them together. To begin with, it will
be convenient to define $L = \frac{2(BC)^2}{AC}$. The first relation
is

\be \lll{rat1}\frac{L(QR)}{L(Pv)} = \frac{QR}{Pv} = \frac{PE}{PC} =
\frac{AC}{PC} \ee
The first equality is obvious, the second equality is because
triangles $Pxv$ and $PCE$ are similar, and the third is Lemma
\ref{ell1}. The second relation is obvious.

\be \lll{rat2} \frac{L(Pv)}{(Gv)(Pv)} = \frac{L}{Gv} \ee
The next is Lemma \ref{ell3}

\be \lll{rat3}\frac{(Gv)(vP)}{(Qv)^2} = \frac{(PC)^2}{(DC)^2} \ee
Next we have

\be \lll{rat4} \frac{(Qv)^2}{(Qx)^2} = 1 \ee
This is only approximately true, but as $Q$ gets very close to $P$,
$\ang QPS$ approaches $\ang RPS$, and it follows from this that
$\frac{(Qv)^2}{(Qx)^2} \lar 1$. Next, we have

\be \lll{rat5} \frac{(Qx)^2}{(QT)^2} = \frac{(EP)^2}{(PF)^2} =
\frac{(CA)^2}{(PF)^2} = \frac{(CD)^2}{(CB)^2} \ee
The first equality is because triangles $EPF$ and $QTx$ are similar,
the second is Lemma \ref{ell1}, and the third is Lemma \ref{ell2}.
We now form a new equation with the relationships
(\ref{rat1})-(\ref{rat5}). The left side of the equation is formed
by multiplying the leftmost parts of (\ref{rat1})-(\ref{rat5}), and
the right side is formed by multiplying the rightmost parts of
(\ref{rat1})-(\ref{rat5}). Much cancelation occurs on the left, and
we get

\be \lll{} \frac{L(QR)}{(QT)^2} =
\frac{L(AC)(PC)^2(CD)^2}{(PC)(Gv)(CD)^2(CB)^2} \ee
Plugging in $L = \frac{2(BC)^2}{AC}$, this is

\be \lll{} \frac{L(QR)}{(QT)^2} =
\frac{2(CB)^2(PC)^2(CD)^2}{(PC)(Gv)(CD)^2(CB)^2} =
\frac{2(PC)}{(Gv)} \ee
As $Q$ approaches $P$, $Gv \lar 2PC$, thus $\frac{2(PC)}{(Gv)} \lar
1$, so we may take

\be \lll{} \frac{L(QR)}{(QT)^2} = 1 \ee
Multiply the equation $L(QR) = (QT)^2$ by $\frac{SP^2}{QR}$ to get

\be \lll{}L(SP)^2 = \frac{(SP)^2(QT)^2}{QR} \ee
Since $L$ is a constant depending on the ellipse, we may write this
as

\be \lll{}(SP)^2 \sim \frac{(SP)^2(QT)^2}{QR} \ee
Now, suppose $Q$ is the point that occurs in the orbit some time $t$
later than $P$. For $Q$ close to $P$ the area of sector $SQP$ is
very close to the area of triangle $SQP$, and thus by the corollary
to Theorem \ref{area} we have $t \sim (QT)(SP)$. Also, by Theorem
\ref{accel} $QR$ is approximately $(1/2)at^2$. Thus, we get

\be \lll{}(SP)^2 \sim \frac{t^2}{at^2}= \frac{1}{a} \ee
Since $SP = R$, the distance from $S$ to $O$, we see $a \sim 1/R^2$, and we are done. \qed

This section should have convinced the reader that a force directed
at the sun satisfying an inverse square law is a likely culprit for
moving the planets in their orbits.

\section{Proof that conical orbits result from an inverse square law}

We now will assume that we have an inverse square law force acting on the planets, and deduce consequences.

\begin{theorem} \lll{unique2} Suppose that an object $O$ moves along a curve
and undergoes an acceleration $a$ towards an unmoving object $S$
which depends only on $r$, the distance between $O$ and $S$. Suppose
that at some time $O$ is at a point $P_1$, and that at some later
time $O$ is at another point $P_2$. For any point $r$ on $SP_1$,
draw the line from $r$ perpendicular to $SP_1$ with length equal to
the acceleration $a(r)$ at distance $r$. The points $a(r)$ trace out
a curve. Let $r_1 = P_1$, and let $r_2$ be the point on $SP_1$ so
that $SP_2 = Sr_2$. Let $v(P)$ be the velocity of $O$ at point $P$.
Then $v(P_2)^2 - v(P_1)^2 = 2A$, where $A$ is the area determined by
the curve $a(r)$ and the line $SP$, between the points $r_1$ and
$r_2$.
\end{theorem}

\includegraphics[width=100mm, height=80mm]{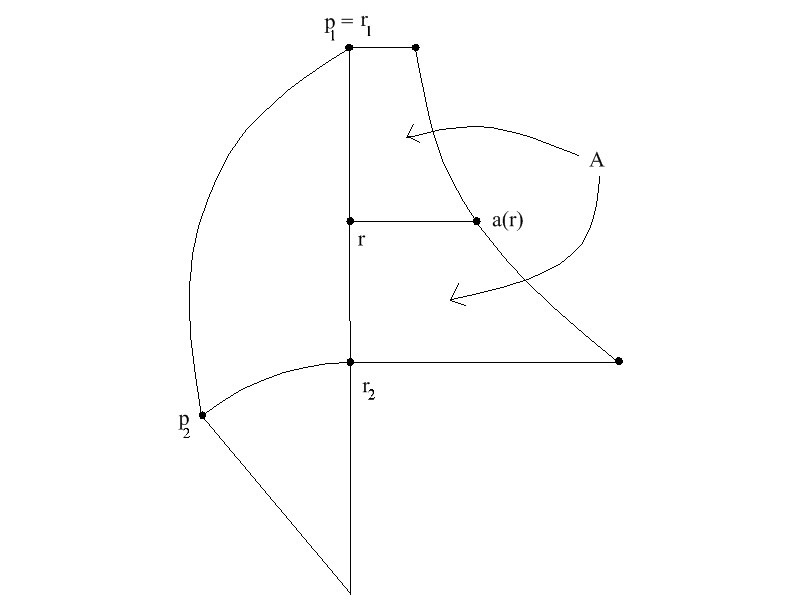}

{\bf Proof:} Divide the interval $[r_2,r_1]$ into $N$ equal parts,
$[x_N,x_{N-1}], [x_{N-1},x_{N-2}],$ $..., [x_2,x_1],[x_1,x_0]$,
where $x_N = r_2$ and $x_0 = r_1$. The following gives an example,
with $N=4$.

\includegraphics[width=100mm, height=80mm]{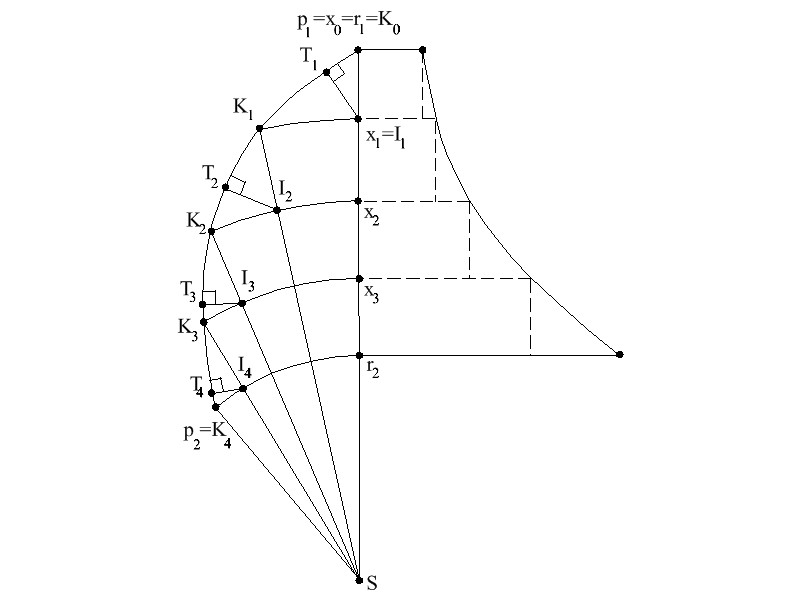}

Let $K_0 = P_1$, $K_N = P_2$, and $K_i$ be the point on the curve
$P_1 P_2$ the same distance from $S$ as $x_i$. Draw the circular arc
between $K_i$ and $x_i$ with center at $S$, and let $I_i$ be the
point on this arc which also lies on $SK_{i-1}$. When $N$ is very
large, the segments $K_{i-1} K_i$ will be approximated as straight
lines. Let $T_i$ be the point on $K_{i-1} K_i$ so that $I_i T_i$ is
perpendicular to $K_{i-1} K_i$. We will approximate by assuming that
the acceleration towards $S$ on $O$ is $a(x_{i-1})$ when $O$ is in
the interval $K_i K_{i-1}$. This acceleration is along the line
$SK_{i-1}$. Resolve the acceleration into perpendicular components,
along $K_{i-1}T_i$ and $T_i I_i$. The acceleration along $T_i I_i$
does not change the velocity, only the direction. Thus, the change
in velocity comes from the acceleration along $K_{i-1}T_i$. Since
$\triangle K_{i-1} K_i I_i \sim \triangle K_{i-1} I_i T_i$,

\be \frac{K_{i-1}T_i}{K_{i-1}I_i} = \frac{K_{i-1}I_i}{K_{i-1}K_i} =
\frac{x_{i-1}x_i}{K_{i-1}K_i} \ee
Thus, the acceleration along $K_i, K_{i-1}$ is $a(x_{i-1})\Big(
\frac{K_{i-1}T_i}{K_{i-1}I_i}\Big) = a(x_{i-1})\Big(
\frac{x_{i-1}x_i}{K_{i-1}K_i} \Big)$

If $t_i$ is the amount of time it takes $O$ to travel from $K_{i-1}$
to $K_i$, then $v(K_i) - v(K_{i-1}) = a(x_{i-1})\Big(
\frac{x_{i-1}x_i}{K_{i-1}K_i} \Big)t_i$. $t_i$ is given by the
distance between $K_{i-1}$ and $K_i$ divided by the velocity over
this interval, which is approximately $v(K_{i-1})$. We see that

\be \lll{}v(K_i) - v(K_{i-1}) = a(x_{i-1})\Big(
\frac{x_{i-1}x_i}{K_{i-1}K_i} \Big)\frac{K_{i-1}K_i}{v(K_{i-1})} \ee
Thus

\be \lll{}v(K_{i-1})(v(K_i) - v(K_{i-1})) = a(x_i)(x_{i-1}x_i) \ee
The right side is equal to the area of one of the rectangles in the
picture above. We see that if we add this expression for all $i$,
the right hand side becomes approximately equal to $A$. Let's
multiply by 2 just for good measure to get

\be \lll{} 2 \sum_{i=1}^N v(K_{i-1})(v(K_i) - v(K_{i-1})) = 2A \ee
Now for a bit of trickery. If we choose $N$ to be very large, then
$v(K_i)$ and $v(K_{i-1})$ will be very close to each other, so we
may write

\bea \lll{}  && 2 \sum_{i=1}^N v(K_{i-1})(v(K_i) - v(K_{i-1}))
\\ && \hh =
\sum_{i=1}^N \Big( v(K_i)(v(K_i) - v(K_{i-1})) + v(K_{i-1})(v(K_i) -
v(K_{i-1})) \Big) \nn \\ \nn && \hh = \sum_{i=1}^N (v(K_i)^2 - v(K_{i-1})^2)
\eea

This last sum is a telescoping sum:

\bea \lll{} && (v(K_1)^2 - v(K_{0})^2) + (v(K_2)^2 - v(K_{1})^2) + \\ \nn && \hh
\ldots + (v(K_{N-1})^2 - v(K_{N-2})^2) + (v(K_N)^2 - v(K_{N-1})^2)
\eea
The
total is $v(K_N)^2 - v(K_0)^2 = v(P_2)^2 - v(P_1)^2$. Therefore,

\be \lll{} v(r_2)^2 - v(r_1)^2 = 2A \ee
which is what we set out to prove. We made many approximations
above, but as $N$ becomes very large the approximations become more
and more accurate, so that in the limit we get equality. \qed

In light of the work we did in the previous section, we need to be
able to calculate the area under the curve given by $a(r) =
\frac{k}{r^2}$, with $k$ a constant.

\begin{proposition} \lll{quad2}Given the setup in the previous theorem, let
$a(r)$ be given by $\frac{m}{r^2}$ for some constant $m>0$. Then $A
= m(\frac{1}{r_2} - \frac{1}{r_1})$.
\end{proposition}

{\bf Proof:} Fundamental theorem of calculus. See also the note at the end of the paper. \qed

Combining Theorem \ref{unique2} and Proposition \ref{quad2} shows
that, if we start an object in motion at a distance $r_o$ from $S$
with velocity $v(r_o)$, then for any other distance $r$ from $S$
that the object attains we have $v(r)^2 - \frac{2m}{r} = v(r_o)^2 -
\frac{2m}{r_o}$. Let's isolate this as a lemma for future reference.

\begin{lemma} \lll{pot} If an object $O$ is in motion about a motionless object $S$
subject only to an acceleration of $\frac{m}{r^2}$ towards $S$, then

\be \lll{}v(r)^2 - \frac{2m}{r} = C \ee
where $C$ is the constant $v(r_o)^2 - \frac{2m}{r_o}$.
\end{lemma}
Corollary 1 from the previous section implies that once an object
begins orbiting $S$ in, say, a counterclockwise direction, it
continues to orbit counterclockwise. That is, it can not reverse
itself at some point and orbit clockwise. We'll assume from now on
that any object $O$ in orbit around $S$ is orbiting
counterclockwise. Let $\aaa$ be the angle between $SO$ and the
tangent to the curve that the object traces; this gives two choices
for $\aaa$, but we will see below that it doesn't matter which we choose, since the important quantity will be the square of the sine of the angle.

\includegraphics[width=100mm, height=80mm]{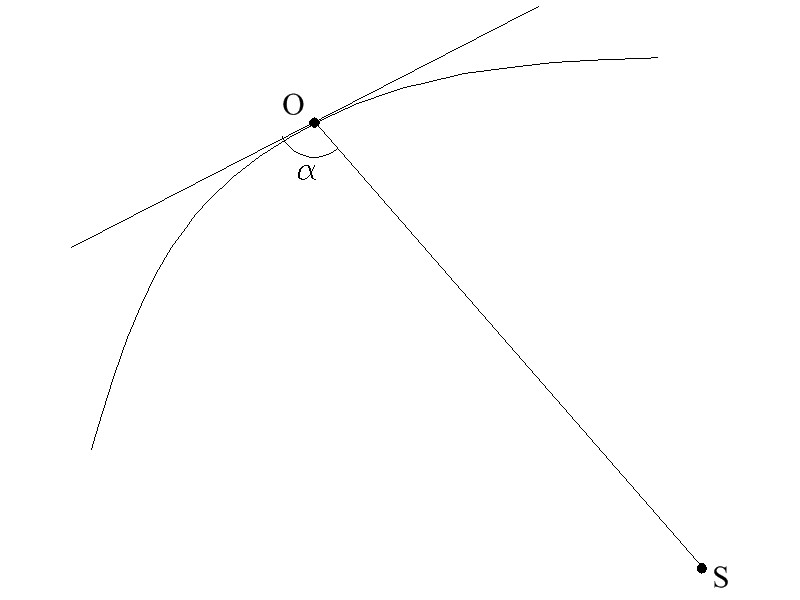}

\begin{lemma} \lll{ai} For an object $O$ in orbit about a motionless object
$S$ with the same assumptions as in Lemma \ref{pot}, we have

\be \lll{}r^2 v(r)^2 \sin^2 \aaa = Q\ee
where $Q$ is a constant.
\end{lemma}

{\bf Proof:} Let $O'$ be a point on the orbit of $O$ close to $O$.
Complete the triangle $OO'S$, and let $\aaa '$ be the angle between
$SO$ and $OO'$, $\bar{\aaa '}$ be the angle between $SO'$ and $O'O$.
Let $t$ be the amount of time $O$ takes to travel to $O'$.

\includegraphics[width=100mm, height=80mm]{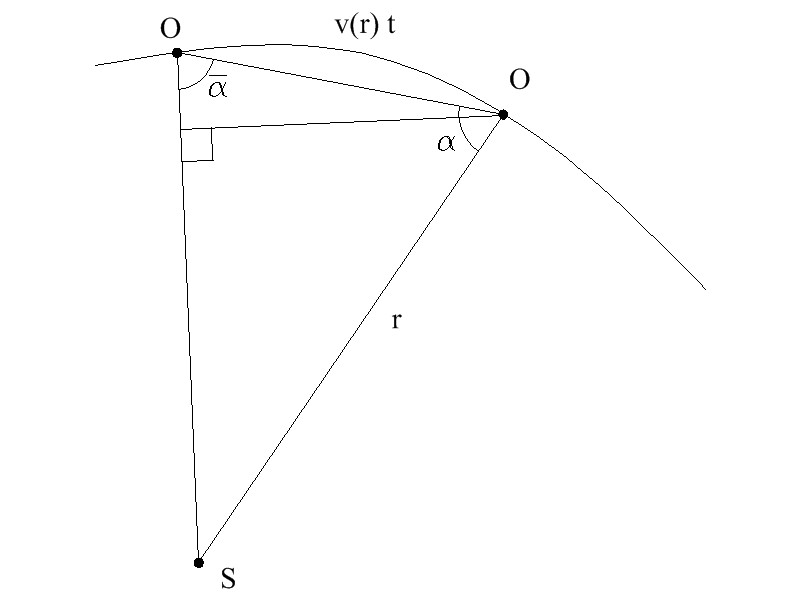}

By Corollary 1 of the previous section, the area of sector $SO'O$ is
$Wt$, for some constant $W$. Triangle $SO'O$ has about the same area
as sector $SO'O$, and $OO'$ is about $v(r)t$, so we approximate

\be \lll{} \frac{1}{2}(SO')(v(r)t)\sin \aaa' = Wt \ee
hence

\be \lll{pen} \frac{1}{2}(SO')v(r)\sin \aaa' = W \ee
As we let $O'$ go to $O$, the approximation becomes equality, $SO'
\lar r$, and $\aaa' \lar \aaa$. Furthermore, $\angle O'SO \lar 0$,
so since $\angle O'SO + \aaa' + \bar{\aaa '} = 180^\circ$,
$\hat{\aaa '} \lar 180^\circ - \aaa$. Substituting these limits into
(\ref{pen}) and using the fact that $\sin (180^\circ - \aaa) = \sin
\aaa$, we get

\be \lll{pen} \frac{1}{2}r v(r)\sin \aaa = W \ee
Squaring this equation and letting $Q = 4W^2$ completes the proof.
\qed

These lemmas combine to give us

\begin{theorem} \lll{jw} An object $O$ in motion about an object $S$ which is
constantly subject to an acceleration of $\frac{m}{r^2}$ in the
direction of $S$ satisfies

\be \lll{oasis} \csc^2 \aaa = \frac{C}{Q}r^2 + \frac{2m}{Q} r \ee
Where $C=, Q,$ and $m$ are constants given earlier in this section,
and $\aaa = \aaa(r)$ is the angle between $SO$ and the tangent to
the orbit at $O$.
\end{theorem}

{\bf Remark:} So that we have the constants all in one place, $m = \frac{a}{r^2}$ which is assumed to be constant, $C= v_o^2 - \frac{2m}{r_o}$, and $Q=r_o^2 v_o^2 \sin^2 \aaa_o$, where $r_o, v_o$ and $\aaa_o$ are the initial distance of $O$ from $S$, velocity of $O$, and angle that $O$ is traveling from the radial line to $S$.
\vski

{\bf Proof:} Rewrite this as $v(r)^2 = \frac{Q}{r^2 \sin^2 \aaa}$
and rewrite the conclusion of Lemma \ref{pot} as $v(r)^2 = C +
\frac{2m}{r}$. Combining these equations gives

\be \lll{} \frac{Q}{r^2 \sin^2 \aaa} = C + \frac{2m}{r} \ee
Multiplying both sides by $\frac{r^2}{Q}$ gives

\be \lll{sec2} \csc^2 \aaa = \frac{C}{Q}r^2 + \frac{2m}{Q} r \ee

\qed

Recall from section 1 that conics satisfy

\be \lll{sec3}\csc ^2 \aaa = \frac{(e^2 - 1)}{e^2(ao)^2}r^2 +
\frac{2}{e(ao)}r \ee
We're clearly getting close, as (\ref{sec2}) and (\ref{sec3}) look
to be pretty much the same equation. There are a few more technical
details to deal with before we can conclude that objects must move
in conic sections under the effect of an inverse square law. The
following proposition is the first step.

\begin{proposition} \lll{circuni} Suppose that an object $O$ is orbiting an object $S$ under
an acceleration of $\frac{m}{r^2}$, and suppose that the orbit of
$O$ contains a circular arc centered at $S$. Then $O$ must move in
that same circle for all times in the past and future.
\end{proposition}

{\bf Proof:} Let $OP$ be the circular arc in the orbit, and complete
the circle around $S$. Let $d$ be the point opposite $O$, let $r_o =
SO$, let $b$ be a point on $OP$ close to $O$, and let $c$ be the
point on the tangent at $O$ so that $bc$ is perpendicular to said
tangent. Suppose that it takes time $t$ for $O$ to move to $b$. We
resolve the motion $Ob$ into components $Oc$ and $cb$. $Oc$ is due
to the velocity at $O$, and thus is equal to $v(r_o)t$. $cb$ is due
to the acceleration towards $S$ at $O$, and is therefore parallel to
$SO$, and is equal to $\frac{m}{2r_o^2}t^2$ by Theorem \ref{accel}.

\includegraphics[width=100mm, height=80mm]{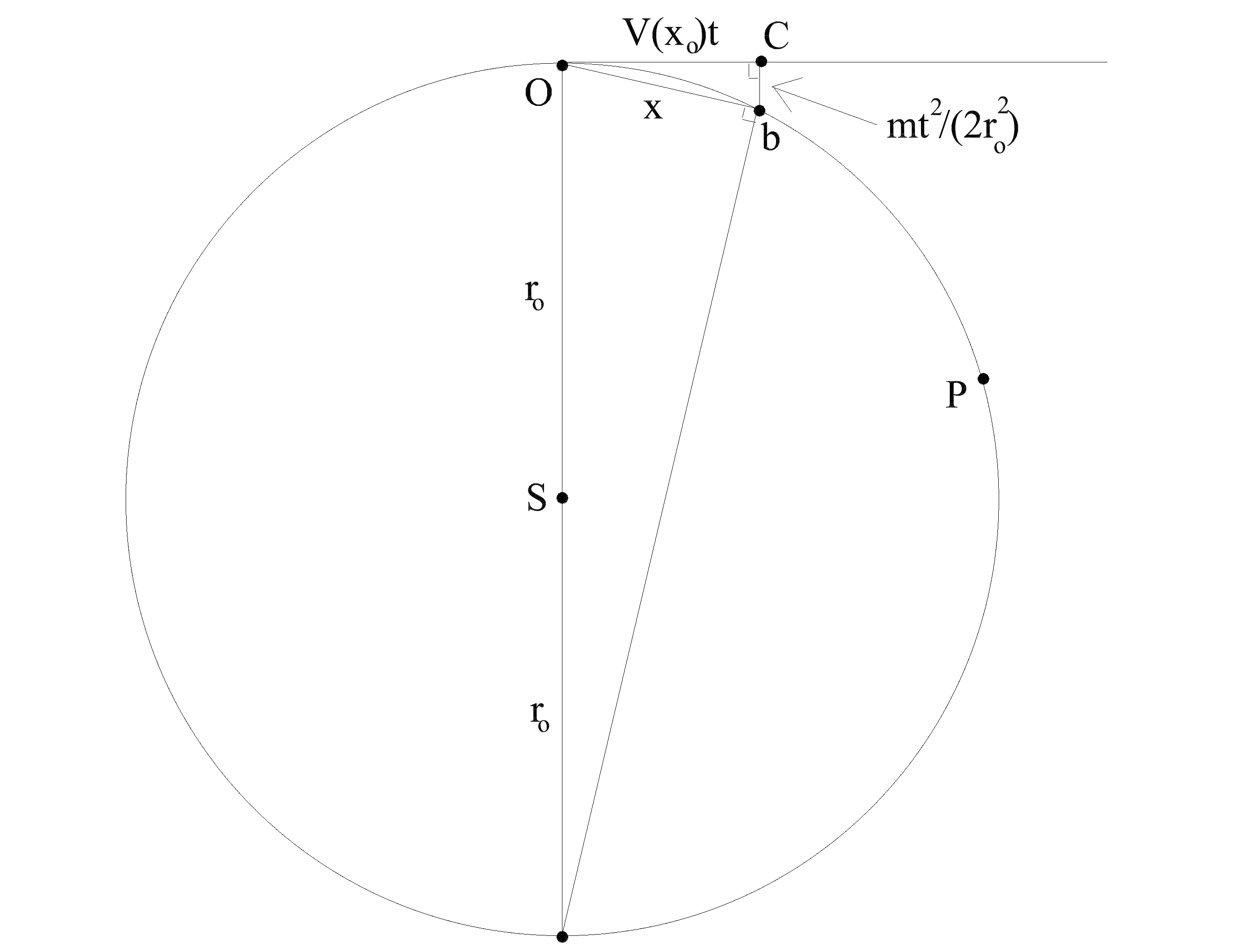}

$\ang Odb$ and $\ang bOc$ both cover the arc $Ob$ and are therefore
equal. Thus, $\tri Odb \sim \tri bOc$, and we have

\bea \lll{smoo} \frac{x}{2r_o} = \frac{\frac{m}{2r_o^2}t^2}{x} \\
\nn \Longrightarrow x^2 = \frac{m}{r_o}t^2 \eea
As $b \lar O$, $\ang bOc \lar 0$, so that $\frac{x}{Oc} \lar 1$. We
can therefore replace $x$ by $v(r_o)t$ in (\ref{smoo}) to get

\be \lll{} v(r_o)^2 t^2 = \frac{m}{r_o}t^2 \Longrightarrow v(r_o)^2
= \frac{m}{r_o} \ee
Now, recall that the $C$ in Lemma \ref{pot} and Theorem \ref{jw} is
given by

\be \lll{} C = v(r_o)^2 - \frac{2m}{r_o} \ee
We see that in our case, $C = \frac{-m}{r_o}$. Now, when $O$ is
moving in a circular arc $\aaa = 90^\circ$, and $\sin 90^\circ = 1$,
so the $Q$ in Lemma \ref{ai} and Theorem \ref{jw} is given by

\be \lll{} Q = v(r_o)^2 r_o^2 \ee
We know that $v(r_o)^2 = \frac{m}{r_o}$, so $Q = mr_o$. Thus, by
Theorem \ref{jw},

\be \lll{} \csc^2 \aaa = \frac{-m/r_o}{mr_o}r^2 + \frac{2m}{mr_o}r =
\frac{-1}{r_o^2}r^2 + \frac{2}{r_o}r = 1 - (\frac{r}{r_o} - 1)^2 \ee
for all $r$ throughout the orbit of $O$. But $\csc ^2 \aaa \geq 1$
for all $\aaa$, and $1 - (\frac{r}{r_o} - 1)^2 \leq 1$ for all $r$.
We see that $\csc^2 \aaa = 1 - (\frac{r}{r_o} - 1)^2 = 1$ for all
$r$ that $O$ can attain. This can only be the case if $r=r_o$, so we
conclude that $r=r_o$ is the only possible distance between $O$ and
$S$. In other words, $O$ moves in a circle of radius $r_o$ forever,
and must have at all times since the object was put in motion. \qed

Now we find the conic that $O$ must travel upon when put in motion.

\begin{proposition}\lll{conicuni} If $O$ is placed into orbit around $S$, then there is
exactly one conic that $O$ can travel along.
\end{proposition}

{\bf Proof:} By Theorem \ref{jw}, the relationship

\be \lll{fire} \csc^2 \aaa = \frac{C}{Q}r^2 + \frac{2m}{Q} r \ee
persists throughout the orbit of $O$. If there is only one possible
$r$ which can satisfy this, then we are in the case covered by the
previous proposition and $O$ travels in a circle, which is a conic.
If $O$ does not travel in a circle, however, we have to find a
unique conic which satisfies (\ref{fire}). In the notation of
Theorem 2 of the first section, a conic(non-circular) is uniquely
determined by the eccentricity $e$ and length $ao$, and satisfies

\be \lll{jc}\csc ^2 \aaa = \frac{(e^2 - 1)}{e^2(ao)^2}r^2 +
\frac{2}{e(ao)}r \ee
Equating coefficients in (\ref{fire}) and (\ref{jc}), we get

\bea \lll{ftf} \frac{e^2 - 1}{e^2(ao)^2} = \frac{C}{Q} \\ \nn
\frac{2}{e(ao)} = \frac{2m}{Q} \eea
The second equation implies

\be \lll{}\frac{1}{e^2(ao^2)} = \frac{m^2}{Q^2} ,\ee

and plugging this into the first gives

\bea \lll{} (e^2-1)\frac{m^2}{Q^2} = \frac{C}{Q}  \\ \nn
\Longrightarrow e^2 = \frac{QC + m^2}{m^2}\eea
Thus,

\be e =\frac{\sqrt{QC + m^2}}{m}. \ee
Given this, the second equation in (\ref{ftf}) implies that

\be ao = \frac{Q}{\sqrt{QC + m^2}}. \ee
So we see that the conic is uniquely determined, i.e. there is only
one conic that satisfies (\ref{fire}). The reader may notice a
possible problem, however. How do we know that $QC + m^2 \geq 0$, so
that we may in fact take the square root? Reexamining the
definitions of the constants, we see that $Q > 0$, but that $C$ can
be any real number. Thus, for arbitrary $Q, C,$ and $m$ it can
easily happen that $QC + m^2 < 0$. What saves us, though, is that
$Q,C,$ and $m$ are not arbitrary. Recall that $C = v(r_o)^2 -
\frac{2m}{r_o}$, and that $Q = v(r_o)^2 r_o^2 \sin^2 \aaa$. If $C
\geq 0$ then we have no problems, so let us assume that $C < 0$.
Then

\bea \nn \lll{} &&  QC + m^2 = \Big( v(r_o)^2 - \frac{2m}{r_o} \Big)
(v(r_o)^2 r_o^2 \sin^2 \aaa) + m^2 \\ && \hskip1.8cm \nn \geq \Big(
v(r_o)^2 - \frac{2m}{r_o} \Big) (v(r_o)^2 r_o^2) + m^2 \\ \nn &&
\hskip1.8cm = r_o^2v(r_o)^4 - 2mr_ov(r_o)^2 + m^2 \\ \nn &&
\hskip1.8cm = (r_ov(r_o)^2 - m)^2 \geq 0 \eea

If $QC + m^2 = 0$ it can be seen easily that $O$ travels in a
circle. Otherwise, $QC + m^2 > 0$, and this entire construction
works to generate a unique conic upon which $O$ can move. \qed

We're almost there. We still have to prove that an object cannot
travel in an orbit that is not a conic. That is, we need to prove
that there is no other curve that satisfies (\ref{oasis}). Try as I might, I couldn't find a geometrical argument for this that varies notably from the standard proof that two functions with the same derivative and same value at a point coincide. I will therefore leave the proof of the following theorem to the reader. If the reader runs into trouble, they might find the transformation $(x,y) \lar (e^x \cos y, e^x \sin y)$ useful, together with the calculus theorem alluded to earlier in this paragraph.

\bt \lll{unique}
Suppose that an object $O$ is placed in orbit around $S$ subject to an equation of the form
\be \label{}
\csc ^2 \aaa = \phi (r)
\ee
\noi Then there is at most one possible path that the object can move along which does not contain circular arcs centered at $S$.
\et

At long last, combining Theorem \ref{jw}, Proposition \ref{circuni}, Proposition \ref{conicuni}, and Theorem \ref{unique}, we obtain

\begin{theorem} \lll{kep1} An object $O$ subject to a force with an inverse
square law directed at an unmoving object $S$ will move along the
path of a conic section.
\end{theorem}

Let's take stock of where we are as it relates to the orbits of the
planets. We have assumed the existence of an acceleration upon any
object in the solar system that is directed at the sun, and which
satisfies an inverse square law. We have proved
that planets and other objects must move along conic
sections(Theorem \ref{kep1}). This gives Kepler's first law, that
planets move in ellipses(if they moved in parabolas or hyperbolas
they would fly out of the solar system, and we wouldn't think of
them as planets). Kepler's second law was proved at the beginning of
the previous section, and in fact would hold for any force directed
at the sun. All that remains is Kepler's third law, which is a snap
compared to the first law. First, one last lemma about ellipses.

\begin{lemma}Let $E$ be an ellipse with focus $a$ and directrix $L$.
Let $o$ be the point on $L$ such that $ao$ is perpendicular to $L$,
and let $e$ be the eccentricity of $E$. Then the area of $E$ is

\be \lll{} \frac{\pi (ao)^2 e^2}{(1-e^2)^{3/2}} \ee

\end{lemma}

{\bf Proof:} Let $X$ and $Y$ denote the major and minor axes of $E$,
respectively. Recall that $E$ can be thought of as the projection of
a circle of radius $Y$, and that projections preserve the ratio of
areas. Let us inscribe $E$ in a rectangle $S'$, and inscribe a
circle $C$ of radius $Y$ in a square $S$.

\hspace{.7cm} \includegraphics[width=100mm, height=80mm]{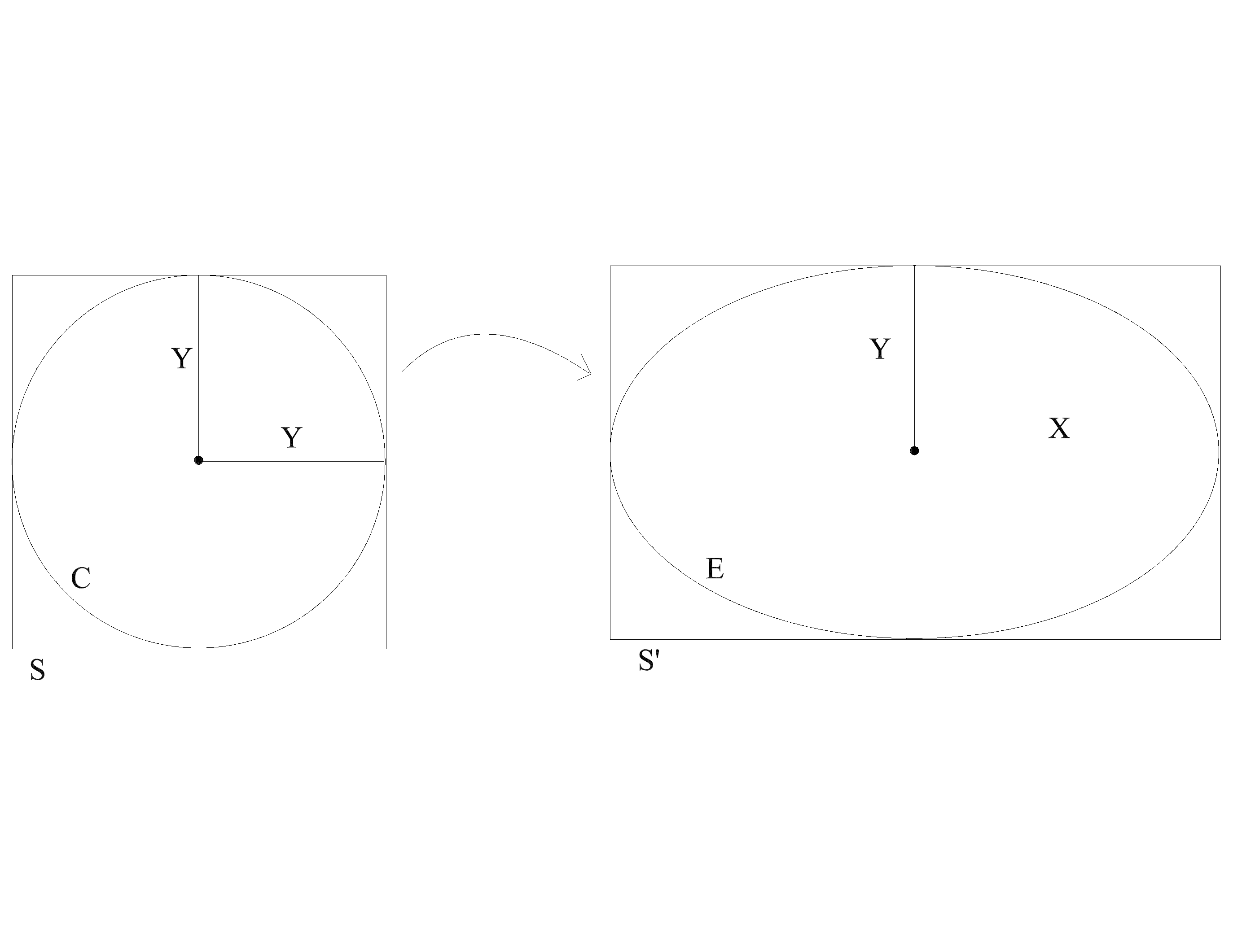}

We know that

\be \lll{}\frac{Area(C)}{Area(S)} = \frac{Area(E)}{Area(S')} \ee
Since $Area(C) = \pi Y^2$, $Area(S) = 4Y^2$, and $Area(S') = 4XY$,
we see that $Area(E) = \pi XY$. Now, by the last proposition in the
section on conics, we have

\be \lll{} X = \frac{(ao)e}{1-e^2} \ee

\be \lll{} Y = \frac{(ao)e}{\sqrt{1-e^2}} \ee
Plugging these identities into $Area(E) = \pi XY$ gives the result.
\qed

That exponent of $3/2$ in the denominator in this lemma sure is
suspicious, isn't it? Now for Kepler's third law.

\begin{theorem} Let an object $O$ orbit a fixed point $S$ in an ellipse $E$ subject
only to a force directed towards $S$ which satisfies an inverse
square law. Let $X$ be the major axis of the ellipse. Then $X^{3/2}
\sim T$, where $T$ is the length of time for $O$ to revolve once
around $S$.
\end{theorem}

{\bf Proof:} Suppose the object is set in motion at time $0$ with
initial velocity $v_o$, radius from $S$ $r_o$, and with angle
between tangent and radius $\aaa_o$. If the planet moves a very
short amount of time $t$ to point $P$ it will sweep out an area that
is very close to a triangle, as below.

\includegraphics[width=100mm, height=80mm]{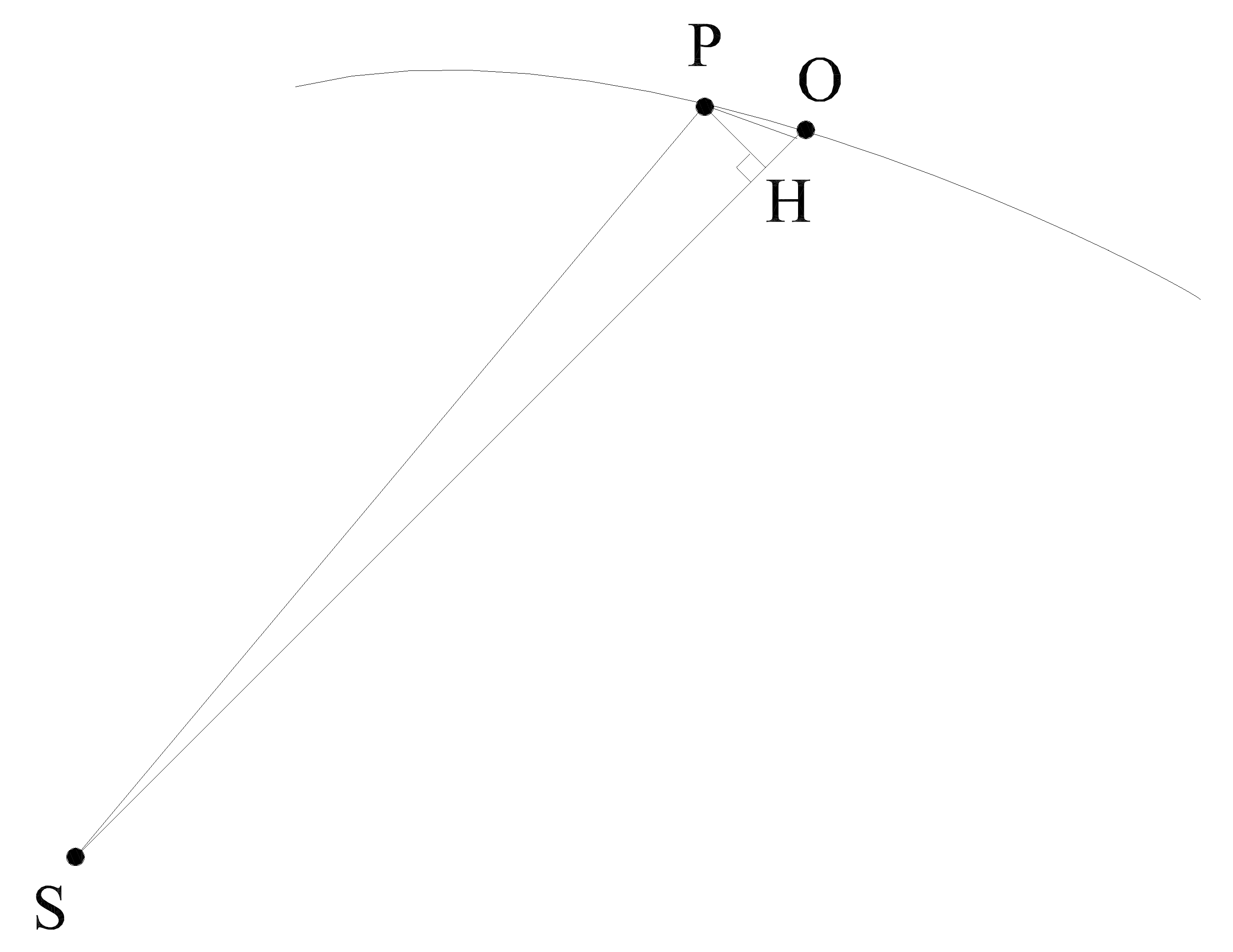}

Since $t$ is so small, $\ang SOP$ is approximately equal to $\aaa$,
and $OP$ is approximately $v_o t$. $SO = r_o$, so the area of
triangle $SOP$ is roughly $tr_o v_o \sin \aaa$. We will take this as
the approximation to the area of sector $SOP$. Since equal areas are
swept out in equal times, this relationship persists throughout the
duration of the orbit. That is, if the planet travels for a length
of time $T$, the radius to $S$ sweeps out the area $Tr_o v_o \sin
\aaa$. Thus, to find the amount of time in one revolution of $O$
about $S$ we may set

\be \lll{jag} Tr_o v_o \sin \aaa = Area(E) \ee
and solve for $T$. From the previous lemma,

\be \lll{} Area(E) = \frac{\pi (ao)^2 e^2}{(1-e^2)^{3/2}} \ee
The proof of Lemma \ref{conicuni} shows that

\be e =\frac{\sqrt{QC + m^2}}{m}. \ee

\be ao = \frac{Q}{\sqrt{QC + m^2}}. \ee
with constants $Q$, $C$, and $m$ as defined earlier in the
section($Q$ defined in Lemma \ref{ai}, $C$ and $m$ in Lemma
\ref{pot}). Thus,

\be \lll{} Area(E) = \frac{Q^2/m^2}{(-QC/m^2)^{3/2}} =
\frac{m\sqrt{Q}}{(-C)^{3/2}} = \frac{mr_o v_o \sin
\aaa_o}{(-C)^{3/2}}\ee
Recall that $O$ moves in an ellipse only when $C < 0$, so that we
may safely raise $(-C)$ to a non-integer power. In light of
(\ref{jag}), we have

\be \lll{ac} T = \frac{m}{(-C)^{3/2}} \ee
We also know that

\be \lll{dd} X = \frac{(ao)e}{1-e^2} = \frac{Q/m}{1-(QC + m^2)/m^2}
= \frac{1}{m}\Big( \frac{1}{-C} \Big) \ee
Comparing (\ref{ac}) and (\ref{dd}) shows that, indeed, $X^{3/2}
\sim T$. \qed

\section{An interesting problem}

{\bf Problem:} {\it Suppose that an object $O$ is placed in orbit around $S$ at a distance $r_o$, an initial velocity $v_o$, and an initial angle $\aaa_o$(assumed not equal to $0^\circ$ or $180^\circ$) to the radial line from $S$. Suppose that $O$ is always subject to an acceleration of $\frac{m}{r^2}$ towards $S$. Determine which of the conic sections $O$ will travel along, determine the closest distance $O$ will attain from $S$, and in the case where $O$ travels in an ellipse determine the maximal distance $O$ attains from $S$.}
\vski

\noi {\bf Remark:} In the case of the parabola and hyperbola, we may need to run time backwards to achieve the minimum, as the object may be placed in motion moving away from $S$.

\vski

\noi {\bf Solution:} We know from Theorem \ref{jw} that, with $C=v_o^2 - \frac{2m}{r_o}$ and $Q=r_o^2 v_o^2 \sin^2 \aaa_o$,
\be \label{}
\csc^2 \aaa = \frac{C}{Q}r^2 + \frac{2m}{Q}r
\ee

From Theorem \ref{gook} we know that this represents a parabola if $v_o^2 - \frac{2m}{r_o}=0$, an ellipse if $v_o^2 - \frac{2m}{r_o}<0$, and a hyperbola if $v_o^2 - \frac{2m}{r_o}>0$. That answers the first part of the problem. To deal with the rest, suppose first that we are in the case of a parabola. Then

\be \label{}
\csc^2 \aaa = \frac{2m}{Q}r
\ee

Since $\csc^2 \aaa \geq 1$, the minimum that $r$ can be is $\frac{Q}{2m}=\frac{r_o^2 v_o^2 \sin^2 \aaa_o}{2m}$. Now suppose that $O$ moves in an ellipse or hyperbola. Again the extremal values of $r$ correspond to $\csc^2 \aaa = 1$, which by the quadratic formula happens when

\be \label{}
r=\frac{-m\pm \sqrt{m^2+CQ}}{C}
\ee

In the proof of Proposition \ref{conicuni} it was shown that $m^2+CQ \geq 0$, so this equation makes sense. When $C>0$ this gives one positive value for $r$, corresponding to the nearest point to $S$ on the hyperbola, and when $C<0$ this gives two positive values, corresponding to the maximal and minimal points to $S$ on the ellipse. Plugging the values in for $C$ and $Q$ gives

\be \label{}
\frac{-m+\sqrt{m^2+(v_o^2 - \frac{2m}{r_o})(r_o^2 v_o^2 \sin^2 \aaa_o)}}{v_o^2 - \frac{2m}{r_o}}
\ee

as the minimum for both the hyperbola and the ellipse, and

\be \label{}
\frac{-m-\sqrt{m^2+(v_o^2 - \frac{2m}{r_o})(r_o^2 v_o^2 \sin^2 \aaa_o)}}{v_o^2 - \frac{2m}{r_o}}
\ee

as the maximum for the ellipse. \qed

The complexity of these answers indicates that a solution by different methods is likely to be quite involved.

\section{Notes, references, and further reading}

1. The prevalence of $r^2$ and $\sin ^2 \aaa$ in the formulas in this paper remind me a bit of {\it rational trigonometry}, as propounded by Norman Wildberger. Essentially this is trigonometry with the fundamental concepts being the squares of lengths and squares of sines of angles. Perhaps many of these theorems could be reworked and would have nicer proofs and statements in that framework. I haven't worked on it myself, but an interested reader might want to consider it. I'm not sure how something like Theorem \ref{gook} would fit in, given the presence of a linear term in $r$.

\vski

\noi 2. The proof of Proposition \ref{han} is based on a technique that, to my knowledge, was discovered by Japanese mathematicians a few centuries ago. I learned of it from \cite{japtemp}, which is highly recommended.

\vski

\noi 3. The Fundamental Theorem of Calculus was invoked only once, in the proof of Proposition \ref{quad2}. But this could have been avoided if one is in a truly classical frame of mind, by a simple argument which is similar to the proof of Theorem \ref{accel}.

\vski

\noi 4. The May, 1994 issue of The College Mathematics Journal contains a very interesting discussion on the question of whether Newton proved the theorem that an inverse square law implies conic section orbits. See also \cite{newtorb}.

\section{Acknowledgements}

I'd like to thank George Markowsky, Linda Markowsky, and Fred Gardiner for their helpful comments and encouragement.

\end{document}